# MULTIVARIATE BAYESIAN FUNCTION ESTIMATION[1]

By Jean-François Angers and Peter T. Kim

*Université de Montréal and University of Guelph*

Bayesian methods are developed for the multivariate nonparametric regression problem where the domain is taken to be a compact Riemannian manifold. In terms of the latter, the underlying geometry of the manifold induces certain symmetries on the multivariate nonparametric regression function. The Bayesian approach then allows one to incorporate hierarchical Bayesian methods directly into the spectral structure, thus providing a symmetry-adaptive multivariate Bayesian function estimator. One can also diffuse away some prior information in which the limiting case is a smoothing spline on the manifold. This, together with the result that the smoothing spline solution obtains the minimax rate of convergence in the multivariate nonparametric regression problem, provides good frequentist properties for the Bayes estimators. An application to astronomy is included.

**1. Introduction.** This paper develops Bayesian function estimation for the multivariate nonparametric regression problem where the domain is taken to be a compact Riemannian manifold. The approach is to incorporate multivariate hierarchical Bayesian methods into the spectral structure of the Riemannian manifold, allowing one to explicitly capture any prior information about possible invariance or symmetry in the data. In particular, a symmetry-adaptive multivariate Bayes estimator is proposed. This approach is a very natural way of modeling symmetries and presents a superior alternative to a frequentist approach. We now discuss this below.

The sample space is usually taken as Euclidean. However, there has been interest in non-Euclidean sample spaces, with the main example being the unit sphere in various dimensions. Without going into a discussion of the history of directional statistics, for that one can consult a very modern account

Received August 2002; revised December 2004.

[1]Supported in part by NSERC Grants R 5297 and OGP 46204.

of the subject in [32], let us provide some details as to the general interest in the subject of a Riemannian sample space, that is, the sample space being a Riemannian manifold, along with the statistical study of certain symmetries inherent in the manifold. In some applications these symmetries often arise due to certain physical constraints imposed by laws of motion.

The greatest number of investigations into the statistical study of symmetries on a Riemannian sample space involves testing. An early paper by Beran [4] investigates testing for uniformity on compact homogeneous spaces, followed by generalizations to compact Riemannian manifolds in [16] and investigation into the two-sample problem on Riemannian manifolds by Wellner [48]. Extensions to specific but more general manifolds have been documented in [8], Chapter 6, which include the case of Grassmann and Stiefel manifolds. Recently, Chikuse and Jupp [9] consider tests of uniformity on shape space. Now one can consider testing for uniformity on a Riemannian manifold as the ultimate form of symmetry. However, "partial" symmetries, such as rotational invariance around a particular axis on a unit sphere, can be of even greater interest in certain physical applications. Starting within the framework of Giné [16], Jupp and Spurr [23] examine testing for symmetries that are not necessarily the full symmetry of uniformity, but only partial symmetries. In more technical terms, which will be made precise below, testing for uniformity can be associated with studying invariance with respect to the full group of isometries on the manifold, whereas partial symmetries can be associated with studying invariance with respect to proper subgroups of the full isometry group.

For multivariate function estimation on Riemannian manifolds, there are some investigations where the primary approach is frequentist. There are a number of works on multivariate function estimation on a unit sphere, see [5, 18, 19, 25, 26, 28, 41, 44, 45, 46]; on Stiefel manifolds, see [29] and [8], Chapter 10; on Lie groups, see [24, 27]; and, on general compact manifolds, see [20, 27]. Other, more exotic analyses have been done in [46], Chapter 2, where multivariate functions are estimated on submanifolds of the tangent bundle manifold of a sphere. There is also some work involving fitting smooth curves on a manifold; see [17, 21, 36]. In some remarkable engineering applications to polymers and robotics, Chirikjian and Kyatkin ([10], Chapters 12 and 17) examine multivariate function estimation on the Euclidean motion group. In computer vision and pattern representation, the sample space is taken to be a certain space of complex matrices, with a pattern being defined through symmetries on this manifold; see [42, 43]. Thus, one can see that there is interest in multivariate function estimation on sample spaces beyond the Euclidean space, and together with current computing capabilities, it is foreseeable that the demand for statistical techniques that go beyond the traditional Euclidean sample space will increase.

The statistical interest in testing for symmetry, partial or full, and in frequentist function estimation on Riemannian manifolds is, nevertheless, disjoint, although it is conceivable that there would be advantages to combining them if possible. Indeed, a frequentist approach would involve one of the following options: assume symmetry (whatever form) is present in the data; ignore any symmetry altogether; or initially test for any symmetries using the methods established in [23], followed by function estimation depending on the outcome of the test. One can see that the above three frequentist approaches have certain difficulties. The first approach runs into the difficulty that symmetry may not really be present, while the second approach has the difficulty that symmetry may be present. The third approach has difficulties in interpreting mean squared error calculations since one is conditioning on that part of the sample space that rejects (or accepts) the symmetry in question. Thus, one can see there are some shortcomings in the frequentist approach. Alternatively, a Bayesian approach to multivariate function estimation on manifolds appears to be a superior solution in that the Bayesian approach allows one to exploit any invariance or symmetry in the data directly, by eliciting very specific prior information on the possible symmetries in question through the spectral structure on the Riemannian manifold. It is in this way that we mean (mentioned in the opening paragraph) that the Bayesian approach is natural for this class of problems and is the subject of this paper.

We now provide a summary of what is to come. In Section 2 notation and some geometric preliminaries are provided. As well, some explicit descriptions of the manifolds discussed above are presented. In Section 3 we initiate a frequentist approach to the multivariate nonparametric regression function estimation problem. This defines the minimization problem in a reproducing kernel Hilbert space; see [46], Chapters 1 and 2. We show that a unique solution defines what may be termed a spline on a manifold which, in turn, confirms a conjecture raised by Wahba [45]. It is shown that the spline solution attains the minimax rate of convergence. As a precursor to the next section, we generalize a result which shows that the spline solution has Bayesian connections when diffuse priors are used. In Section 4 we formally embark upon the task at hand by incorporating initial prior information into the model. We treat the model as involving both a symmetric and nonsymmetric part, assume normality on the first stage finite-dimensional priors and treat the infinite-dimensional part of the model as a nuisance parameter. We can then control the amount of symmetry by controlling the variance terms and employ hyperpriors to deal with the prior parameters. This leads to symmetry-adaptive hierarchical Bayesian function estimators. Bayes factors are subsequently used to determine the truncation level. As for dealing with the nuisance parameters, we can diffuse some of them away and in so doing, we can obtain as limits the smoothing spline solutions on manifolds. This suggests that the hierarchical Bayes estimator has good frequentist properties, which are discussed in Section 5. In Section 6 we go through a detailed analysis for the 2-sphere, as well as present some numerical work on long period cometary orbits. It is here that we see the benefits of incorporating prior symmetries into the model. This extends an earlier study of this data set; see [23]. Section 7 contains the proofs and an appendix is included for further needed technical details.

**2. Notation.** Let $\mathbb{M}$ be a compact connected orientable Riemannian manifold. Consider the Riemannian structure $\{g(p) : p \in \mathbb{M}\}$ and let $dx$ be the normalized volume element of $\mathbb{M}$ associated with this structure. For each fixed $p \in \mathbb{M}$, we can associate with $g(p)$ a matrix $(g_{ij}(p))$ called the metric tensor. We will, in addition, assume that the manifold is without boundary, although we could generalize the following arguments to certain boundary conditions.

Let $C^\infty(\mathbb{M})$ be the space of real-valued infinitely differentiable continuous functions on $\mathbb{M}$. Denote by

$$\Delta = -\frac{1}{\sqrt{|g(p)|}} \sum_{j,k} \partial_j (g^{jk}(p) \sqrt{|g(p)|} \partial_k)$$

the Laplace–Beltrami operator on $\mathbb{M}$, where $\partial_j$ denotes the partial derivative with respect to the $j$th component, $(g^{ij}(p))$ is the inverse of $(g_{ij}(p))$ the metric tensor and $|g(p)|$ is the determinant of the matrix $(g_{ij}(p))$. We note that $\Delta$ is an elliptic self-adjoint second-order differential operator on $C^\infty(\mathbb{M})$ for which the eigenfunctions of $\Delta$ are a complete orthonormal basis for $\mathcal{L}^2(\mathbb{M})$, the space of real-valued square integrable functions.

Let $\lambda$ be an eigenvalue of $\Delta$. The collection of all eigenvalues for a given $\mathbb{M}$ is countably infinite, hence, letting $\mathbb{N}_0 = \{0, 1, 2, \ldots\}$, we can enumerate the eigenvalues by $\lambda_k \geq 0$, $k \in \mathbb{N}_0$, with no upper bound. Furthermore, we will use the convention that $\lambda_0 = 0$ and that $\lambda_k \leq \lambda_{k+1}$ for $k \in \mathbb{N}_0$. For each $k$, let $\phi_{kj}$ be an eigenfunction so that $\Delta \phi_{kj} = \lambda_k \phi_{kj}$, for $j = 1, \ldots, L_k$, and denote by $\mathcal{E}_k = \mathrm{sp}\{\phi_{kj} : j = 1, \ldots, L_k\}$, where $\mathrm{sp}(\cdot)$ stands for the span of the object in question. Then $\dim \mathcal{E}_k = L_k < \infty$, $k \in \mathbb{N}_0$, where $\dim(\cdot)$ denotes the dimension of the object in question.

Let $\phi_k = (\phi_{k1}, \ldots, \phi_{k \dim \mathcal{E}_k})'$, where superscript "$\prime$" denotes transpose and let $\langle \cdot, \cdot \rangle_k$ denote the dot product on $\mathbb{R}^{\dim \mathcal{E}_k}$, with $\|\cdot\|_k$ the induced norm, $k \in \mathbb{N}_0$. For $h \in \mathcal{L}^2(\mathbb{M})$, the eigenfunction expansion will be defined by

$$(2.1) \qquad h = \sum_{k=0}^{\infty} \langle \hat{h}_k, \phi_k \rangle_k \qquad \text{where } \hat{h}_k = \int_{\mathbb{M}} h \phi_k,$$

for $k \in \mathbb{N}_0$, where integration over $\mathbb{M}$ is defined piecewise using the usual partition of unity argument; see (A.1) in the Appendix.

We can consider subspaces of $\mathcal{L}^2(\mathbb{M})$ in the following way. First, on the space $C^\infty(\mathbb{M})$ of infinitely continuous differentiable functions on $\mathbb{M}$, consider the Sobolev norm $\|\cdot\|_{H_s}$ of order $s$ defined accordingly. For any function $h = \sum_k \langle \hat{h}_k, \phi_k \rangle_k$, let

$$\|h\|_{H_s}^2 = \sum_{\lambda_k > 0} \lambda_k^s \|\hat{h}_k\|_k^2. \tag{2.2}$$

One can verify that (2.2) is indeed a norm. Denote by $H_s(\mathbb{M}) \subset \mathcal{L}^2(\mathbb{M})$ the (vector-space) completion of $C^\infty(\mathbb{M})$ with respect to (2.2). This will be called the Sobolev space of order $s > \dim \mathbb{M}/2$. In addition, we will also consider Sobolev ellipsoids defined by

$$H_s(\mathbb{M}, M) = \left\{ u \in H_s(\mathbb{M}) : \sum_k \lambda_k^s \|\hat{u}_k\|_k^2 \leq M \right\},$$

for $M > 0$ and $s > \dim \mathbb{M}/2$.

Often $\mathbb{M}$ is equipped with certain symmetries. A Riemannian manifold is homogeneous if its group of isometries $\mathbb{G}$ acts transitively on $\mathbb{M}$, where by the latter we mean that for every $p, q \in \mathbb{M}$ there exists an $h \in \mathbb{G}$ such that $p = hq$, where multiplication denotes the group action $\mathbb{G} \times \mathbb{M} \to \mathbb{M}$. For every $p \in \mathbb{M}$, let $\mathbb{G}^p = \{h \in \mathbb{G} : hp = p\}$ denote the isotropy subgroup of $p$. It is well known that if $\mathbb{M}$ is a homogeneous compact connected Riemannian manifold, then, for every $p \in \mathbb{M}$, $\mathbb{G}^p$ is a closed subgroup of $\mathbb{G}$ and there exists a diffeomorphism between the quotient space $\mathbb{G}/\mathbb{G}^p$ and $\mathbb{M}$. A differentiable function $f : \mathbb{M} \to \mathbb{R}$ is called a zonal function with respect to the action of the isotropy subgroup $\mathbb{G}^p$, $p \in \mathbb{M}$, if it is constant on the isotropy subgroup $\mathbb{G}^p$.

New manifolds can be created by taking products of existing manifolds, as well as by taking quotients as in the case of a homogeneous space. For concreteness, let us consider specific examples that will be used for illustration throughout the paper, as well as those mentioned in Section 1.

EXAMPLE 2.1 (Sphere). The sphere $S^{p-1} \subset \mathbb{R}^p$ is the set of unit vectors in $p$-dimensional Euclidean space. In the case where $p = 3$, we note that any point on $S^2$ can almost surely be represented by

$$\omega = (\cos \varphi \sin \vartheta, \sin \varphi \sin \vartheta, \cos \vartheta)', \tag{2.3}$$

where $\varphi \in [0, 2\pi)$, $\vartheta \in [0, \pi)$ and superscript "$\prime$" denotes transpose.

EXAMPLE 2.2 (Orthogonal and special orthogonal group). The orthogonal group $O(p)$ consists of the space of $p \times p$ real orthogonal matrices. However, this group is not connected. The connected component consisting of those real orthogonal matrices having determinant equal to unity, $SO(p)$, is called the special orthogonal group. Again, in the case of $p=3$, $SO(3)$ can be represented in the following way. Let

$$u(\varphi) = \begin{pmatrix} \cos\varphi & -\sin\varphi & 0 \\ \sin\varphi & \cos\varphi & 0 \\ 0 & 0 & 1 \end{pmatrix}, \qquad a(\vartheta) = \begin{pmatrix} \cos\vartheta & 0 & \sin\vartheta \\ 0 & 1 & 0 \\ -\sin\vartheta & 0 & \cos\vartheta \end{pmatrix},$$

where $\varphi \in [0, 2\pi)$, $\vartheta \in [0, \pi)$. The well-known Euler angle decomposition says any element of $SO(3)$ can almost surely be uniquely written as

$$g = u(\varphi_1)a(\vartheta)u(\varphi_2),$$

where $\varphi_1 \in [0, 2\pi)$, $\varphi_2 \in [0, 2\pi)$, $\vartheta \in [0, \pi)$.

New manifolds can be created from products and quotients of these examples. In fact, it turns out that $SO(3)$ is the transitive group of isometries on $S^2$. Furthermore, the subgroup

$$SO(2) = \{u(\varphi) : \varphi \in [0, 2\pi)\}$$

of $SO(3)$ is the isotropy subgroup of $(0,0,1)' \in S^2$. Thus, we can identify $S^2$ with the quotient space $SO(3)/SO(2)$; hence, $S^2$ is a homogeneous space. Throughout the paper we will use the sphere $S^2$, as well as its transitive group of isometries, $SO(3)$, for illustrative purposes leading up to the application in Section 6.

Some of the other manifolds previously mentioned include the Stiefel manifold, $V_k(\mathbb{R}^p) = O(p)/O(p-k)$, the Grassmann manifold, $G_k(\mathbb{R}^p) = O(p)/(O(p) \times O(p-k))$, and shape space, $S^{p(k-1)-1}/SO(p)$. Of the more exotic constructions, the collection of all tangent spaces on a manifold is called the tangent bundle and is itself a manifold. The Euclidean motion group is defined to be $SO(p) \ltimes \mathbb{R}^p$, where $\ltimes$ denotes a semi-direct product. In computer vision, the sample space is taken to be $SL(2, \mathbb{C})/O(2)$, where $SL(2, \mathbb{C})$ denotes the space of $2 \times 2$ complex matrices of determinant 1. We note that the last two examples are that of noncompact manifolds.

As for orthonormal bases, we have the following example.

EXAMPLE 2.3 (Spherical harmonics). Let

$$(2.4) \quad Y_{kq}(\omega) = \begin{cases} \sqrt{2}\sqrt{\dfrac{(2k+1)(k-q)!}{4\pi(k+q)!}} P_q^k(\cos\vartheta)\cos(q\varphi), & q = 1,\ldots,k, \\ \sqrt{\dfrac{(2k+1)}{4\pi}} P_0^k(\cos\vartheta), & q = 0, \\ \sqrt{2}\sqrt{\dfrac{(2k+1)(k-|q|)!}{4\pi(k+|q|)!}} P_{|q|}^k(\cos\vartheta)\sin(|q|\varphi), & \\ & q = -1,\ldots,-k, \end{cases}$$

where $\varphi \in [0, 2\pi)$, $\vartheta \in [0, \pi)$ and $P_q^k$ are the Legendre functions, $-k \leq q \leq k$ and $k \in \mathbb{N}_0$. We note that we can think of (2.4) as the vector entries to the $(2k+1)$-vector

$$Y_k(\omega) = (Y_{kq}(\omega)),$$

$|q| \leq k$ and $k \in \mathbb{N}_0$. In this situation $\{Y_{kq} : |q| \leq k, k \in \mathbb{N}_0\}$ are the eigenfunctions of the Laplace–Beltrami operator on $S^2$ with eigenvalues $\lambda_k = k(k+1)$, $k \in \mathbb{N}_0$, and, hence, form a complete orthonormal basis over $\mathcal{L}^2(S^2)$.

Some further technical properties are provided in the Appendix. In addition, the following asymptotic notation will be used. Let $\{a_n\}$ and $\{b_n\}$ denote two real sequences of numbers. We write $a_n \ll b_n$ to mean $a_n \leq Cb_n$ for some $C > 0$ as $n \to \infty$, the Vinogradov notation. This notation is more convenient than the "big oh" Landau notation since expressions that are within an order of magnitude can be long and, furthermore, can create confusion with the notation for the orthogonal group. We will, however, use the notation $a_n = o(b_n)$ to mean $a_n/b_n \to 0$ as $n \to \infty$. Furthermore, $a_n \asymp b_n$ whenever $a_n \ll b_n$ and $b_n \ll a_n$, and $a_n \sim b_n$ when $a_n/b_n \to 1$ as $n \to \infty$.

**3. Nonparametric regression and splines on manifolds.** Let $f \in \mathcal{L}^2(\mathbb{M})$. Its eigenfunction expansion, as defined in (2.1), is

(3.1) $$f(x) = \sum_{k=0}^{\infty} \langle \gamma_k, \phi_k(x) \rangle_k,$$

for $x \in \mathbb{M}$ and where

(3.2) $$\gamma_k = \int_{\mathbb{M}} f(x) \phi_k(x) \, dx, \qquad \phi_k \in \mathcal{E}_k,$$

for $k = 0, 1, \ldots$. If we observe (3.1) at the points $x_1, \ldots, x_n \in \mathbb{M}$, then our observations are

(3.3) $$y_i = f(x_i) + \varepsilon_i \qquad \text{for } i = 1, \ldots, n,$$

where we assume $\varepsilon = (\varepsilon_1, \ldots, \varepsilon_n)' \sim N(0, \sigma^2 I)$, $\sigma^2 > 0$, and we are interested in estimating $f$, a real-valued function on $\mathbb{M}$. We note that the Fourier coefficients in (3.1) and (3.2) are denoted by $\gamma_k$ and not $\hat{f}_k$, $k \in \mathbb{N}_0$, as was done in Section 2. The purpose for this departure is due to the fact that the Bayesian framework will later treat the coefficients as random quantities.

For any fixed value of $K > 0$, called the truncation level, (3.1) can be written as

(3.4) $$f(x) = f_K(x) + \eta(x),$$

where

$$(3.5) \qquad f_K(x) = \sum_{k=0}^{K} \langle \gamma_k, \phi_k(x) \rangle_k \quad \text{and} \quad \eta(x) = \sum_{k=K+1}^{\infty} \langle \gamma_k, \phi_k(x) \rangle_k.$$

Then we can write the regression problem (3.3) as

$$(3.6) \qquad y = \Phi\gamma + \eta + \varepsilon,$$

where $y = (y_1, \ldots, y_n)'$, $\eta = (\eta(x_1), \ldots, \eta(x_n))'$, $\gamma = (\gamma_1, \ldots, \gamma_\kappa)'$, $\kappa = \sum_{k=0}^{K} \dim \mathcal{E}_k$, $\varepsilon = (\varepsilon_1, \ldots, \varepsilon_n)'$ and $\Phi = (\phi_k(x_i))$ for $k = 0, 1, \ldots, K$, $i = 1, \ldots, n$.

3.1. *Splines on manifolds.* We will need the following notation. Let $x_1, \ldots, x_n$, $x \in \mathbb{M}$, and define

$$Q^\lambda(x_{i_1}, x_{i_2}) = \sum_{k \geq K+1} \lambda_k^{-s} \langle \phi_k(x_{i_1}), \phi_k(x_{i_2}) \rangle_k.$$

Define the $n \times n$ matrix

$$Q^\lambda_{n,\xi} = [Q^\lambda(x_{i_1}, x_{i_2})] + n\xi I_n,$$

where $x_{i_1}, x_{i_2} \in \mathbb{M}$, $i_1, i_2 = 1, \ldots, n$, $\xi \geq 0$ and $I_n$ is the $n \times n$ identity matrix. Furthermore, define the $\kappa \times 1$ vector $\phi(x) = (\phi_k(x))$, and define the $n \times 1$ vector $q(x) = [Q^\lambda(x_1, x), \ldots, Q^\lambda(x_n, x)]'$. The following generalizes earlier multivariate spline smoothing methods of Wahba [45], Cox [11] and Taijeron, Gibson and Chandler [41].

THEOREM 3.1. *Let $\mathbb{M}$ be a compact connected orientable Riemannian manifold. Assume $x_1, \ldots, x_n$ are distinct points on $\mathbb{M}$, $x \in \mathbb{M}$, and consider the following smoothing problem:*

$$\min_{u \in H_s(\mathbb{M})} \frac{1}{n} \sum_{i=1}^{n} (u(x_i) - y_i)^2 + \xi \int_{\mathbb{M}} |\Delta^{s/2} u(x)|^2 \, dx,$$

*where $\xi > 0$, $s > \dim(\mathbb{M})/2$ and $y_i \in \mathbb{R}$ for $i = 1, \ldots, n$. Define*

$$f_\xi^n(x) = \phi(x)'d + q(x)'c,$$

*where the $n \times 1$ vector $c$ and the $\kappa \times 1$ vector $d$ are defined by*

$$c = [Q^\lambda_{n,\xi}]^{-1}(I_n - \Phi(\Phi'[Q^\lambda_{n,\xi}]^{-1}\Phi)^{-1}\Phi'[Q^\lambda_{n,\xi}]^{-1})y,$$

$$d = (\Phi'[Q^\lambda_{n,\xi}]^{-1}\Phi)^{-1}\Phi'[Q^\lambda_{n,\xi}]^{-1}y,$$

*with $y = (y_1, \ldots, y_n)'$. Then $f_\xi^n(x)$ for $\xi > 0$ is the unique solution to the smoothing problem.*

REMARK 3.2. The practical choice of the smoothing parameter $\xi > 0$ can be determined by using generalized cross-validation similar to that outlined for the Euclidean case in, for example, [46], Chapter 4.

REMARK 3.3. As stated in Section 1, one can always create new manifolds by taking products of existing ones. This leads to tensor splines as the multivariate functions to be estimated over product manifolds. Such is the approach taken in [31], which is a special case of the more general tensor spline construction; see [46], Chapter 10. In particular, one can freely construct any number of products of compact manifolds since the individual reproducing kernels can be amalgamated into one. Thus, for example, one can take products of various spheres, rotation matrices, Stiefel and Grassmann manifolds.

3.2. *Minimaxity.* In terms of assessing the spline estimator of Theorem 3.1, we are led to ask about the type of frequentist properties they possess. Details of minimaxity for spline estimators have been investigated in the univariate case by Speckman [40], along with an extension to a multivariate framework by Cox [11, 12].

We first state the following sharp lower bound result which is directly related to Theorem 1 of [35] and is stated as Theorem 2.1 of [14].

THEOREM 3.4 (Pinsker and Efromovich). *Let $\mathbb{M}$ be a compact connected orientable manifold without boundary. Then*

$$\inf_{\tilde{f}} \sup_{f \in H_s(\mathbb{M}, M)} \mathbb{E}\|\tilde{f} - f\|^2$$
$$\geq (M\mathcal{W}^{2s/\dim \mathbb{M}})^{\dim \mathbb{M}/(2s+\dim \mathbb{M})} \wp n^{-2s/(2s+\dim \mathbb{M})}(1 + o(1))$$

*as $n \to \infty$, where the infimum is taken over all estimators, $s > \dim \mathbb{M}/2$,*

$$\mathcal{W} = \mathcal{W}(\mathbb{M}) = \frac{\mathrm{vol}\,\mathbb{M}}{(2\sqrt{\pi})^{\dim \mathbb{M}} \Gamma(1 + \dim \mathbb{M}/2)},$$
$$\wp = \wp(s, \dim \mathbb{M})$$
$$= \left(\frac{2s}{2s + 2\dim \mathbb{M}}\right)^{2s/(2s+\dim \mathbb{M})} \left(\frac{2s + \dim \mathbb{M}}{\dim \mathbb{M}}\right)^{\dim \mathbb{M}/(2s+\dim \mathbb{M})},$$

vol$(\cdot)$ *denotes volume and $\Gamma(\cdot)$ denotes the gamma function.*

REMARK 3.5. We emphasize that this lower bound is over all estimators. The constant $\mathcal{W}$ is a geometric invariant associated with asymptotic calculations performed by Hermann Weyl; see (A.2) in the Appendix. The constant $\wp$ is associated with asymptotic calculations performed by Pinsker [35].

The existence of the spline solution only requires the design points $x_1, \ldots, x_n \in \mathbb{M}$ to be distinct. More, however, will be required in order to calculate the integrated mean squared error. In the Euclidean univariate case, Speckman ([40], page 972) states the necessary condition on the design needed to achieve the minimax lower bound. For the multivariate case, the assumptions are listed in [11], pages 791 and 792. In both of the above univariate and multivariate cases, the domain in question is embedded in a Euclidean space of the same dimension. This, however, is not possible for a general compact manifold; hence, the assumptions necessary on the design points $x_1, \ldots, x_n \in \mathbb{M}$ must satisfy local versions of the assumptions of Cox ([11], page 791).

Let us state the four assumptions for the Euclidean case. Assumption 1 is just the statement of the nonparametric regression model (3.3). For $U \subset \mathbb{R}^d$, a bounded open set, consider the points $t = (t_1, \ldots, t_d)', t_k = (t_{k1}, \ldots, t_{kd})' \in U \subset \mathbb{R}^d$, $k = 1, \ldots, n$. Define the empirical distribution function in the usual way,

$$F_n(t) = \sum_{t_{kj} \leq t_j} n^{-1},$$

where summation occurs over every coordinate, $j = 1, \ldots, d$ and $k = 1, \ldots, n$. Let $F(t)$ denote the cumulative distribution function for $t_1, \ldots, t_n$ and define

$$b_n = \sup_t |F_n(t) - F(t)|.$$

Assumption 2 states that the smoothing parameter, $\xi > 0$, must satisfy

$$\xi \in [\xi_n, \Xi_n], \qquad \xi_n \leq \Xi_n, \qquad \lim_{n \to \infty} b_n \xi_n^{-5 \dim \mathbb{M}/(4s)} = \lim_{n \to \infty} \Xi_n = 0.$$

Furthermore, assumption 3 says that the density of $F(t)$ must be bounded away from 0 and infinity, while assumption 4 requires the open set $U \subset \mathbb{R}^d$ to be bounded simply connected with a smooth boundary.

The local version of assumptions 2, 3 and 4 can go as follows. By the definition of a manifold, for every point $p \in \mathbb{M}$, there exist an open set $\mathcal{O}_\alpha \subset \mathbb{M}$, a local diffeomorphism $\psi_\alpha : \mathcal{O}_\alpha \to \psi_\alpha(\mathcal{O}_\alpha) \subset \mathbb{R}^{\dim \mathbb{M}}$ and $p \in \mathcal{O}_\alpha$. The pair $(\mathcal{O}_\alpha, \psi_\alpha)$ is called a chart, and the collection of all charts is called an atlas if it covers $\mathbb{M}$. The atlas is defined in greater detail in the Appendix. The local version of assumptions 2 and 3 of [11], page 791, is, therefore, that these assumptions take place on every chart. Integration on $\mathbb{M}$ requires a partition of unity $\mathcal{P} = \{\delta_\alpha : \alpha \in \mathcal{A}\}$, subordinate to the atlas where $\delta_\alpha : \mathbb{M} \to [0, 1]$ with the support, $\mathrm{supp}\, \delta_\alpha \subset \mathcal{O}_\alpha$, $\alpha \in \mathcal{A}$; see the Appendix. The local version of Cox's assumption 4 is that the boundary of $\mathrm{supp}\, \delta_\alpha$ is smooth. To give these three assumptions a name, we will say that $x_1, \ldots, x_n \in \mathbb{M}$ locally satisfies the Cox assumptions. We have the following result.

THEOREM 3.6. *Suppose $x_1, \ldots, x_n \in \mathbb{M}$ locally satisfies the Cox assumptions. If $\xi \asymp n^{-2s/(2s+\dim \mathbb{M})}$, then*

$$\mathbb{E} \|f_\xi^n - f\|^2 \ll n^{-2s/(2s+\dim \mathbb{M})}$$

*as $n \to \infty$, where $f \in H_s(\mathbb{M}, M)$ for some $M > 0$ and $s > (\dim \mathbb{M})(5 \dim \mathbb{M} - 2)/4$.*

Thus, in terms of the rate of convergence, the spline estimator of Theorem 3.1, with stronger order of smoothness, achieves the lower bound. We suspect that the required order of smoothness can substantially be reduced and this will be pursued elsewhere.

3.3. *Splines as Bayes estimators with diffuse priors.* Although smoothing splines are a general computational method for function fitting, there is, however, a Bayesian interpretation. This Bayesian approach to smoothing splines on the unit interval is discussed in [46], Chapter 1. In the following we adapt that approach for $\mathbb{M}$ partly to generalize Theorem 1.5.3 of [46], but mainly because subsequent hierarchical modeling builds from this earlier work.

To generalize Theorem 1.5.3 of [46], we need to consider the concept of a random field $X$ on $\mathbb{M}$. In general, the random field can be expanded in terms of the eigenfunctions so that

$$X(p) = \sum_k \langle Z_k, \phi_k(p) \rangle_k, \qquad p \in \mathbb{M},$$

where $Z_k$ is a sequence of independent, mean zero, $\dim \mathcal{E}_k$-dimensional random vectors, with each coordinate having variance $\sigma_k^2$, $k \in \mathbb{N}_0$. If, in addition, each $Z_k$, $k \in \mathbb{N}_0$, is normally distributed, we say that the process $X$ is Gaussian with covariance kernel

$$EX(p)X(q) = \sum_{k \geq 0} \sigma_k^2 \langle \phi_k(q), \phi_k(p) \rangle_k,$$

for $p, q \in \mathbb{M}$ and $\sigma_k \geq 0$, $k \in \mathbb{N}_0$. We have the following result.

THEOREM 3.7. *Let $X(x)$ be a mean zero real-valued Gaussian random field defined on $\mathbb{M}$ with covariance kernel*

$$Q^\lambda(x_1, x_2) = \sum_{k=K+1}^{\infty} \lambda_k^{-s} \langle \phi_k(x_1), \phi_k(x_2) \rangle_k,$$

*for $x_1, x_2 \in \mathbb{M}$, for some $K > 0$ and $s > \dim \mathbb{M}/2$. Consider*

$$f(x) = \sum_{k=0}^{K} \langle \gamma_k, \phi_k(x) \rangle_k + \tau X(x) \qquad \text{for } x \in \mathbb{M} \text{ and } \tau > 0,$$

*and suppose we observe $f$ at the points $x_1, \ldots, x_n \in \mathbb{M}$. Let our observations be*

$$y_i = f(x_i) + \varepsilon_i \qquad \text{for } i = 1, \ldots, n,$$

*where $\varepsilon = (\varepsilon_1, \ldots, \varepsilon_n)' \sim N(0, \sigma^2 I)$, $\sigma^2$ known. Suppose $\gamma_k | \tau^2 \sim N(0, \tau^2 \nu I_{\dim \mathcal{E}_k})$, $\tau^2$ and $\nu$ known, for $k = 0, 1, \ldots, K$. Consider the Bayes solution*

$$\tilde{f}_\nu(x) = \mathbb{E}(f(x) | y_1, \ldots, y_n)$$

*and suppose $f_\xi^n(x)$ is the solution to the smoothing problem with $\xi = \sigma^2/(n\tau^2)$. Then for each fixed $x \in \mathbb{M}$,*

$$\lim_{\nu \to \infty} \tilde{f}_\nu(x) = f_\xi^n(x).$$

**4. Symmetry and Bayesian modeling.** A Riemannian manifold can exhibit symmetries which we want to directly capture in the modeling process. Indeed, as in [23], let $\mathbb{G}^0$ be a subgroup of the isometry group of $\mathbb{M}$; see Section 2. We say that $f \in \mathcal{L}^2(\mathbb{M})$ is invariant under the action of $\mathbb{G}^0$ if $f(gx) = f(x)$ for all $g \in \mathbb{G}^0$ and $x \in \mathbb{M}$. In terms of the eigenstructure, for any $k \in \mathbb{N}_0$ the eigenspace $\mathcal{E}_k$ decomposes into two orthogonal subspaces. Denote by $\mathcal{E}_k^0$ the eigenfunctions in $\mathcal{E}_k$ that are invariant with respect to $\mathbb{G}^0$, which will be referred to below as having $\mathbb{G}^0$-invariance and $\mathcal{E}_k^1$, its orthogonal complement in $\mathcal{E}_k$, which will be referred to below as non-$\mathbb{G}^0$-invariance. This allows us to write

$$\mathcal{E}_k = \mathcal{E}_k^0 \oplus \mathcal{E}_k^1, \tag{4.1}$$

and inner products to be written as

$$\langle \cdot, \cdot \rangle_k = \langle \cdot, \cdot \rangle_k^0 + \langle \cdot, \cdot \rangle_k^1$$

for $k \in \mathbb{N}_0$.

Let $f \in \mathcal{L}^2(\mathbb{M})$. Thus, to exhibit this $\mathbb{G}^0$-invariance let us rewrite (3.1) as

$$f(x) = \sum_{k=0}^\infty \{\langle \gamma_k^0, \phi_k(x) \rangle_k^0 + \langle \gamma_k^1, \phi_k(x) \rangle_k^1\}, \tag{4.2}$$

for $x \in \mathbb{M}$ and where

$$\gamma_k^j = \int_\mathbb{M} f(x) \phi_k(x) \, dx, \qquad \phi_k \in \mathcal{E}_k^j,$$

for $j = 0, 1$, $k \in \mathbb{N}_0$.

REMARK 4.1. We would like to remark that the splitting up of the sums in terms of the $\mathbb{G}^0$-invariant part allows us to later incorporate explicit prior assumptions of $\mathbb{G}^0$-invariance. We note that if one assumes that $\mathbb{G}^0$ is the trivial subgroup, then $\mathcal{E}_k^0 = \mathcal{E}_k$, hence $\mathcal{E}_k^1 = \{0\}$ for all $k \in \mathbb{N}_0$.

EXAMPLE 4.2. In terms of $S^2$, we can consider the subgroup to be $SO(2)$ of the full transitive group $SO(3)$. Thus, for a function $f:S^2 \to \mathbb{R}$, $SO(2)$-invariance would mean that $f$ is a zonal function and thus only depends on $\vartheta \in [0, \pi)$. In terms of (4.1), by looking at the definition of (2.4), one can see that, for $\mathcal{L}^2(S^2)$, $SO(2)$-invariance would mean

(4.3) $$\mathcal{E}_k^0 = \{Y_{k0}\} \quad \text{and} \quad \mathcal{E}_k^1 = \{Y_{kq} : 1 \leq |q| \leq k\}$$

for each $k \in \mathbb{N}_0$.

There are two approaches for dealing with the parameters $\gamma_k$ for $k \geq K+1$. One can engage in eliciting very informative prior information in order to estimate the $\eta_i$'s; see [1, 2]. Alternatively, one can adopt an approach wherein $\eta_i$ are combined with the measurement errors $\varepsilon_i$ for $i = 1, \ldots, n$; see [3]. The latter approach is truly Bayesian, but at the inference stage we can treat these $\eta_i$ as nuisance quantities and eliminate them by integrating out (rather than estimating or diffusing) the corresponding parameters. In this section we will deal with the latter approach. The former approach will be discussed in Section 5 since this method of analysis allows one to make direct comparisons with splines.

4.1. *Eliciting prior information.* Our prior belief in the $\mathbb{G}^0$-invariance of (3.1) under the subgroup $\mathbb{G}^0$, which is explicitly invoked in (4.2), can be captured using a mixture normal model, that is,

(4.4) $$\gamma_k^j | \tau^2 \sim p N(0, \tau^2 \Gamma_k^{j0}) + (1-p) N(0, \tau^2 \Gamma_k^{j1}),$$

for some $\tau^2 > 0$, where $\Gamma_k^{jr} = \text{diag}(\beta_{k\ell}^{jr})$, $\ell = 1, \ldots, \dim \mathcal{E}_k$, $k \in \mathbb{N}_0$, $j = 0, 1$, $r = 0, 1$, and $p$ models our prior belief that the $\mathbb{G}^0$-invariance assumption is true.

The $\mathbb{G}^0$-invariance assumption is also taken into account by having smaller variances. Indeed, the components of the variance of $\gamma_k^1$ will be smaller than that of the components of the variance of $\gamma_k^0$. In addition, we would like to have finite variance for the $\eta_i$'s. All of these properties can be obtained by assuming $\beta_{k\ell}^{11} \leq \beta_{k\ell}^{01}$ and $\beta_{k\ell}^{10} \equiv 0$ for $\ell = 1, \ldots, \dim \mathcal{E}_k$, $k = 1, \ldots, K$. Furthermore, assume $\beta_{k\ell}^{j0}, \beta_{k\ell}^{j1} \leq \lambda_k^{-s}$ for $k = K+1, \ldots$, $j = 0, 1$, where the $\lambda_k$'s are the eigenvalues of the Laplace–Beltrami operator $\Delta$ on $\mathbb{M}$ defined in Section 2.

Once a joint prior distribution is specified for $\sigma^2$ and (the hyperparameter) $\tau^2$, the prior model is complete. Note further that, since later we assign a second stage prior on the variance factor $\tau^2$, their marginal prior distribution will no longer be normal, but a heavier tailed distribution ensuring a certain degree of prior robustness to our estimator (cf. [6], Chapter 4).

Consider the $\kappa \times 1$ vector $\gamma = (\gamma_k^j)$ for $k = 0, 1, \ldots, K$ and $j = 0, 1$. The prior specified above indicates that

$$\gamma | \tau^2 \sim p N(0, \tau^2 \Gamma^0) + (1-p) N(0, \tau^2 \Gamma^1),$$

where the $\kappa \times \kappa$ covariance matrices

$$\Gamma^r = \bigoplus_{j,k} \Gamma_k^{jr},$$

where $r = 0, 1$, with the direct sums being taken over $k = 0, 1, \ldots, K$ and $j = 0, 1$.

For the remainder term, write

$$\eta^j(x) = \sum_{k=K+1}^{\infty} \langle \gamma_k^j, \phi_k(x) \rangle_k^j,$$

for $x \in \mathbb{M}$, $j = 0, 1$. Now

$$(\eta_1^j, \ldots, \eta_n^j)' | \tau^2 \sim p N(0, \tau^2 Q_n^{\Omega j 0}) + (1-p) N(0, \tau^2 Q_n^{\Omega j 1}),$$

where $Q_n^{\Omega jr} = (Q^{\Omega jr}(x_{i_1}, x_{i_2}))$, and

(4.5) $$Q^{\Omega jr}(x_{i_1}, x_{i_2}) = \sum_{k=K+1}^{\infty} \langle \phi_k(x_{i_1}), \Omega_k^{jr} \phi_k(x_{i_2}) \rangle_k^r,$$

where $\Omega_k^{jr} = \text{diag}(\beta_{k\ell}^{jr})$, $\ell = 1, \ldots, \dim \mathcal{E}_k^j$, $k = K+1, \ldots$, $j = 0, 1$, $r = 0, 1$, $x_{i_1}, x_{i_2} \in \mathbb{M}$ and $i_1, i_2 = 1, \ldots, n$. We have the following result.

LEMMA 4.3. *Suppose $0 \le \beta_{k\ell}^{jr} \le \lambda_k^{-s}$, for $r = 0, 1$, $j = 0, 1$, $\ell = 1, \ldots, \dim \mathcal{E}_k^j$, $k = K+1, \ldots$, and $s > \dim(\mathbb{M})/2$, where the $\lambda_k$'s are the eigenvalues of the Laplace–Beltrami operator $\Delta$ on $\mathbb{M}$. Then*

$$|\text{Cov}(\eta_{i_1}^j, \eta_{i_2}^j)| \le \tau^2 C(\mathbb{M}, s),$$

*for $i_1, i_2 = 1, \ldots, n$, $j = 0, 1$, where $C(\mathbb{M}, s) < \infty$ is a constant depending only on $\mathbb{M}$ and $s > \dim(\mathbb{M})/2$.*

4.2. *The posterior.* Consider the $n \times \kappa$ design matrix

$$\Phi = (\phi_k^r(x_i)),$$

where $k = 0, 1, \ldots, K$, $r = 0, 1$ and $i = 1, \ldots, n$. Then we obtain the following structure. Given $\gamma$, $\sigma^2$ and $\tau^2$, we have the following linear model for the observations $y = (y_1, \ldots, y_n)'$:

(4.6) $$y = \Phi \gamma + u,$$

where $u \sim N(0, \Sigma)$ with $\Sigma = \sigma^2 I_n + \tau^2 Q_n^\Omega(p)$, where $Q_n^\Omega(p) = p(Q_n^{\Omega 00} \oplus Q_n^{\Omega 10}) + (1-p)(Q_n^{\Omega 01} \oplus Q_n^{\Omega 11})$. This follows from the fact that

$$y|\gamma, \eta, \sigma^2, \tau^2 \sim N(\Phi\gamma, \sigma^2 I_n) \quad \text{and} \quad \eta|\tau^2 \sim N(0, \tau^2 Q_n^\Omega(p)).$$

From (4.6) and using standard hierarchical Bayes techniques (cf. [30]) and matrix identities (cf. [39], page 151), it follows that

$$
\begin{aligned}
(4.7) \quad y|\sigma^2, \tau^2 \sim\ & pN(0, \sigma^2 I_n + \tau^2(\Phi\Gamma^0\Phi' + Q_n^\Omega(p))) \\
& + (1-p)N(0, \sigma^2 I_n + \tau^2(\Phi\Gamma^1\Phi' + Q_n^\Omega(p))),
\end{aligned}
$$

$$(4.8) \quad \gamma|y, \sigma^2, \tau^2 \sim p^* N(A^0 y, B^0) + (1-p^*) N(A^1 y, B^1),$$

where

$$p^* = \frac{pm^0(y)}{pm^0(y) + (1-p)m^1(y)},$$

$$
\begin{aligned}
(4.9) \quad A^r &= \tau^2 \Gamma^r \Phi'(\sigma^2 I_n + \tau^2(\Phi\Gamma^r\Phi' + Q_n^\Omega(p)))^{-1}, \\
B^r &= \tau^2 \Gamma^r - \tau^4 \Gamma^r \Phi'(\sigma^2 I_n + \tau^2(\Phi\Gamma^r\Phi' + Q_n^\Omega(p)))^{-1} \Phi\Gamma^r,
\end{aligned}
$$

where $r = 0, 1$ and $m^0(y)$, $m^1(y)$ denote, respectively, the normal density with mean vector 0 and covariance matrices $\sigma^2 I_n + \tau^2(\Phi\Gamma^0\Phi' + Q_n^\Omega(p))$ and $\sigma^2 I_n + \tau^2(\Phi\Gamma^1\Phi' + Q_n^\Omega(p))$.

At this point (4.8) allows us to produce an estimator of (3.4) once the hyperparameters in $A^0$ and $A^1$ are set. Two possible ways of handling this are the following: first, to use diffuse prior parameters; or, second, treat the current priors as first stage priors and add additional hyperprior assumptions. The first approach produces generalized Bayes estimators. In the following section we will use the hyperprior approach.

4.3. *Hierarchical Bayesian modeling.* In order to proceed to the second stage calculations, some algebraic simplifications are needed (cf. [1]). Spectral decomposition yields $\Phi\Gamma^r\Phi' + Q_n^\Omega(p) = H^r D^r H^{r\prime}$, where $D^r = \text{diag}(d_1^r, d_2^r, \ldots, d_n^r)$ is the matrix of eigenvalues and $H^r$ the orthogonal matrix of eigenvectors for $r = 0, 1$. Thus,

$$
\begin{aligned}
\sigma^2 I_n + \tau^2(\Phi\Gamma^r\Phi' + Q_n^\Omega(p)) &= H^r(\sigma^2 I_n + \tau^2 D^r) H^{r\prime} \\
&= \tau^2 H^r(v I_n + D^r) H^{r\prime},
\end{aligned}
$$

where $r = 0, 1$ and $v = \sigma^2/\tau^2$. Using this spectral decomposition, the marginal density of $y$ given $\tau^2$ and $v$ can be written as

$$m(y|\tau^2, v) = pm^0(y|\tau^2, v) + (1-p)m^1(y|\tau^2, v),$$

where
$$m^r(y|\tau^2 v) = (2\pi\tau^2)^{-n/2}\det(vI_n + D^r)^{-1/2}$$

(4.10)
$$\times \exp\left\{-\frac{1}{2\tau^2}y'H^r(vI_n + D^r)^{-1}H^{r\prime}y\right\}$$

$$= (2\pi\tau^2)^{-n/2}\left\{\prod_{i=1}^{n}(v+d_i^r)^{-1/2}\right\}\exp\left\{-\frac{1}{2\tau^2}\sum_{i=1}^{n}\frac{(w_i^r)^2}{v+d_i^r}\right\},$$

where $w^r = (w_1^r, \ldots, w_n^r)' = H^{r\prime}y$ for $r = 0, 1$.

EXAMPLE 4.4. In the case of $S^2$, using the spherical harmonics (2.4) and invoking $SO(2)$-invariance on $\mathcal{L}^2(S^2)$, we can decompose the eigenspace as in (4.3). Thus, $\dim \mathcal{E}_k^0 = 1$ and $\dim \mathcal{E}_k^1 = 2k$ for $k = K+1, \ldots$. Furthermore, $\kappa = \sum_{k=0}^{K}(\dim \mathcal{E}_k^0 + \dim \mathcal{E}_k^1) = \sum_{k=0}^{K}(2k+1) = (K+1)^2$. This type of symmetry is often observed in directional data.

Fix
$$\nu_{k\ell}^{j1} = [(k+1/2)(k+1)(k+2)(k+3)]^{-1}$$

for all $j = 0, 1$, $|\ell| \leq k$, $k = K+1, \ldots$ and
$$\nu_{k\ell}^{j0} = \begin{cases} \nu_{k\ell}^{j1}, & \text{if } j = 0, \\ 0, & \text{otherwise.} \end{cases}$$

Then
$$Q^{\Omega 01}(\omega_1, \omega_2)$$

(4.11)
$$= (2\pi)^{-1}\left[(\tfrac{1}{2}q_2(\omega_1'\omega_2) - \tfrac{1}{6})\right.$$
$$\left. - \sum_{k=1}^{K}((k+1/2)(k+1)(k+2)(k+3))^{-1}P_k(\omega_1'\omega_2)\right],$$

where $\omega_1, \omega_2 \in S^2$,

(4.12)
$$q_2(w) = \frac{1}{2}\left\{\ln\left(1 + \sqrt{\frac{2}{1-w}}\right)\left[12\left(\frac{1-w}{2}\right)^2 - 4\left(\frac{1-w}{2}\right)\right]\right.$$
$$\left. - 12\left(\frac{1-w}{2}\right)^{3/2} + 6\left(\frac{1-w}{2}\right) + 1\right\}$$

for $|w| \leq 1$ and $P_k$ is the $k$th Legendre polynomial, $k \in \mathbb{N}_0$. Similarly,
$$Q^{\Omega 00}(\omega, \nu) = \sum_{k=K+1}^{\infty}\frac{P_k(\omega_1'\omega_2)}{(k+1/2)(k+1)(k+2)(k+3)}.$$

4.4. *Second stage prior and estimation.* To derive the function estimator, all that is now needed is to eliminate the hyper and nuisance parameters from the first stage posterior distribution, by integrating out these variables with respect to their second stage prior. Since it is well known (cf. [6], Chapter 4) that the final Bayes estimator does not depend crucially on the second and higher stage hyperpriors, these priors can be chosen to simplify computations. Accordingly, the priors on $\tau^2$ and $v$ can be chosen as $\pi_{2,1}(\tau^2) \propto (\tau^2)^{-c}$; see [3]. With this choice of prior on $\tau^2$, the marginal prior on $\gamma$ has the form

$$\pi(\gamma) \propto p\left(\sum_{k=0}^{K}\sum_{j=0}^{1}\sum_{\ell=1}^{\dim \mathcal{E}_k^j} \frac{(\gamma_{k\ell}^j)^2}{\beta_{k\ell}^{j0}}\right)^{-0.5(\kappa+2c-2)}$$
$$+ (1-p)\left(\sum_{k=0}^{K}\sum_{j=0}^{1}\sum_{\ell=1}^{\dim \mathcal{E}_k^j} \frac{(\gamma_{k\ell}^j)^2}{\beta_{k\ell}^{j1}}\right)^{-0.5(\kappa+2c-2)}.$$

This prior density corresponds to the limiting case of a multivariate Student-$t$ density (which has heavier tails than the likelihood function). The prior on $v$ is chosen to be an $F$-distribution with $a$ and $b$ degrees of freedom satisfying the following conditions:

- the prior variance of $v$ $(= \frac{2b^2(a+b-2)}{a(b-4)(b-2)^2})$ is infinite;
- the Fisher information number $(= \frac{a^2(b+2)(b+6)}{2(a-4)(a+b+2)})$ is minimum;
- the prior mode $(= \frac{b(a-2)}{a(b+2)})$ is greater than 0.

This can be done by choosing $2 < b \leq 4$ and $a = 8(b+2)/(b-2)$. Let $\pi_{22}(v)$ denote the resulting prior density.

Once the second stage priors are specified, using (4.10) and taking the expectation with respect to $\tau^2$, the Bayesian estimator of $\gamma$ under $\mathbb{G}^0$-invariance ($r = 0$) or non-$\mathbb{G}^0$-invariance ($r = 1$) is given by

$$(4.13) \qquad \widetilde{\gamma}^r = \Gamma^r \Phi' H^r \mathbb{E}^r[(vI_n + D^r)^{-1}|y] H^{r'}y,$$

and the expectation is taken with respect to

$$(4.14) \quad \pi_{22}^r(v|y) \propto \frac{v^{a/2-1}}{(b+av)^{(a+b)/2}}\left(\prod_{i=1}^{n}(v+d_i^r)\right)^{-1/2}\left(\sum_{i=1}^{n}\frac{(w_i^r)^2}{v+d_i^r}\right)^{-(n+2c-2)/2},$$

for $r = 0, 1$. Note that in order for $\pi_{22}^r(v|y)$, $r = 0, 1$, to be proper densities, $c$ should be chosen such that $c < b/2$. Hence, under the mixture prior (4.4) and squared error loss, the Bayes estimator for $\gamma$ is given by

$$(4.15) \qquad \widetilde{\gamma} = p^*\widetilde{\gamma}^0 + (1-p^*)\widetilde{\gamma}^1.$$

Again, using (4.10), the posterior expected squared error loss of $\gamma$ can be written as

$$\mathrm{Var}(\gamma|y) = p^* \mathrm{Var}^0(\gamma|y) + (1-p^*)\mathrm{Var}^1(\gamma|y) + 2p^*(1-p^*)(\widetilde{\gamma}^0 - \widetilde{\gamma}^1)(\widetilde{\gamma}^0 - \widetilde{\gamma}^1)',$$

where

$$\begin{aligned}
\mathrm{Var}^r(\gamma|y) &= \frac{1}{n+2c-4}\mathbb{E}^r\left[\sum_{i=1}^n \frac{(w_i^r)^2}{v+d_i} \,\Big|\, y\right] \\
&\quad - \frac{1}{n+2c-4}\Gamma^r \Phi' H^r \mathbb{E}^r\left[\left(\sum_{i=1}^n \frac{(w_i^r)^2}{v+d_i^r}\right)(vI_n + D^r)^{-1}\,\Big|\,y\right] H^{r\prime} \Phi \Gamma^r \\
&\quad + \mathbb{E}^r[\widetilde{\gamma}^r(v)\widetilde{\gamma}^r(v)'|y],
\end{aligned}$$

$r = 0, 1$. Since these expectations involve only one-dimensional integrals, they can be computed easily using one of the several standard techniques, such as Gauss quadrature, Monte Carlo or Laplace approximation.

Finally, the function estimator $\tilde{f}$ of $f$ is

$$\tilde{f}(x) = p^* \tilde{f}^0(x) + (1-p^*)\tilde{f}^1(x),$$

where

$$\tilde{f}^0(x) = \sum_{k=0}^K \langle \widetilde{\gamma}^0, \phi_k(x) \rangle_k^0 \tag{4.16}$$

and

$$\tilde{f}^1(x) = \sum_{k=0}^K \langle \widetilde{\gamma}^1, \phi_k(x) \rangle_k^1, \tag{4.17}$$

for $x \in \mathbb{M}$. Note that equation (4.16) corresponds to the Bayes estimator of $f$ if one believes that the $\mathbb{G}^0$-invariance assumption is true, that is, if $p = 1$, and (4.17) is the Bayes estimator of $f$ under the non-$\mathbb{G}^0$-invariance assumption.

4.5. *Bayes factor and choice of $K$.* We now describe how the optimal level of truncation $K$ is to be determined. As indicated above in (3.4), the choice of $K$ provides a model for the observations through the choice of the corresponding regression function. Denote the maximum truncation level by $K_{\max}$, so that $\sum_{k=0}^{K_{\max}}[\dim \mathcal{E}_k^0 + \dim \mathcal{E}_k^1] < n$ (cf. Section 3).

Let $\mathfrak{M}_K$ denote the model arising from (3.4), (3.5) corresponding to the truncation level $K$. Our task is to pick the best model for the given data from the set of models

$$\mathfrak{M}_K, \qquad K = 1, 2, \ldots, K_{\max}.$$

The well-accepted method (cf. [37], Section 7.2.2) for deciding between two possible models is to compute their associated Bayes factor. As a basis of comparison, the larger model $\mathfrak{M}_{K_{\max}}$ will be used. Hence, we have to compute

$$\mathfrak{B}_K = \frac{m(y|\mathfrak{M}_K)}{m(y|\mathfrak{M}_{K_{\max}})}, \tag{4.18}$$

where $m(y|\mathfrak{M}_K)$ denotes the marginal density of $y$ under the model $\mathfrak{M}_K$, $K = 1, \ldots, K_{\max}$. (Note that $\mathfrak{B}_{K_{\max}} = 1$.) From (4.7) it follows that, under $\mathfrak{M}_K$,

$$y|\sigma^2, \tau^2 \sim pN(0, \sigma^2 I_n + \tau^2(\Phi_K \Gamma_K^0 \Phi_K' + Q_{n,K}^\Omega(p))) \\ + (1-p)N(0, \sigma^2 I_n + \tau^2(\Phi_K \Gamma_K^1 \Phi_K' + Q_{n,K}^\Omega(p))),$$

where we have shown the dependence of $\Phi$, $\Gamma^0$, $\Gamma^1$ and $Q_n^\Omega$ on $K$ explicitly with subscripts. It follows then that

$$m(y|\mathfrak{M}_K) = \int m(y|\mathfrak{M}_K, \sigma^2, \tau^2) \, d\pi(\sigma^2, \tau^2),$$

where $\pi(\sigma^2, \tau^2)$ is the joint prior distribution on $\sigma^2$ and $\tau^2$.

As in the previous section, consider the spectral decomposition of $\Phi_K \Gamma_K^r \Phi_K' + Q_{n,K}^\Omega(p)$, for $r = 0, 1$. Let $D_K^r$ and $H_K^r$ be such that $\Phi_K \Gamma_K^r \Phi_K' + Q_{n,K}^\Omega(p) = H_K^r D_K^r H_K^{r\,\prime}$ for $r = 0, 1$. Also, let $d_{K,i}^r$ be the $i$th diagonal element of $D_K^r$ and let $w_K^r = H_K^r y = (w_{K,1}^r, \ldots, w_{K,n}^r)'$ for $r = 0, 1$. Then, using (4.14), the marginal density of $y$ under $\mathfrak{M}_K$ can be expressed as

$$m(y|\mathfrak{M}_k) = p \int m^0(y|\mathfrak{M}_k, \sigma^2, \tau^2) \, d\pi(\sigma^2, \tau^2) \\ + (1-p) \int m^1(y|\mathfrak{M}_k, \sigma^2, \tau^2) \, d\pi(\sigma^2, \tau^2),$$

where

$$m^r(y|\mathfrak{M}_k) = \int_0^\infty \frac{v^{a/2-1}}{(b+av)^{(a+b)/2}} \\ \times \left(\prod_{i=1}^n (v + d_{K,i}^r)\right)^{-1/2} \left(\sum_{i=1}^n \frac{(w_{K,i}^r)^2}{v + d_{K,i}^r}\right)^{-(n+2c-2)/2} dv,$$

for $r = 0, 1$. Consequently, to choose the best model $\mathfrak{M}_K$, one needs to compute $m(y|\mathfrak{M}_K)$ for $K = 1, \ldots, K_{\max}$. Then the best value of $K$ (equivalently, the best model $\mathfrak{M}_K$) is the one for which $\mathfrak{B}_K$ is maximum.

An alternative to the Bayes factor is to use Schwarz's criterion (cf. [38] and [37], Section 7.2.3) which can be viewed as an approximation to the logarithm of the Bayes factor. Since $\mathfrak{M}_1 \subset \mathfrak{M}_2 \subset \cdots \subset \mathfrak{M}_{K_{\max}}$, Schwarz's criterion can be written as

$$\mathfrak{S}_{ij} = -\log \mathfrak{L}_n + \frac{\pi_j - \pi_i}{2} \log(n), \tag{4.19}$$

where $\mathfrak{L}_n$ denotes the ratio of the likelihood functions under $\mathfrak{M}_i$ and $\mathfrak{M}_j$ evaluated at the maximum likelihood estimator of $\gamma$ under both models and $\pi_i$ ($\pi_j$) corresponds to the number of parameters in model $\mathfrak{M}_i$ ($\mathfrak{M}_j$). The model $\mathfrak{M}_i$ is preferred to the model $\mathfrak{M}_j$ if $\mathfrak{S}_{ij} > 0$.

**5. Frequentist properties of the hierarchical Bayes estimator.** We may wish to model the nuisance parameter directly as in [1, 2]. One of the advantages of this approach is that it allows one to make a direct comparison with the smoothing spline approach of Section 3.1. Although we can approach this by imposing $\mathbb{G}^0$-invariance priors as done previously, in order to ease the notation we will not impose $\mathbb{G}^0$-invariance conditions, or, equivalently, as stated in Remark 4.1, we will assume invariance with respect to the trivial subgroup.

5.1. *Modeling the nuisance parameter directly.* We begin by rewriting (3.6) as

$$y = \Upsilon\psi + \varepsilon,$$

where $\Upsilon = [\Phi, I]$, $\psi = (\gamma', \eta')'$. Following Angers and Delampady [1], we assume a multivariate normal prior

$$\psi \sim N(\psi_0, \tau^2 \Xi), \qquad \Xi = \begin{bmatrix} \Gamma & 0 \\ 0 & Q_n^\Omega \end{bmatrix},$$

where $Q_n^\Omega = (Q^\Omega(x_i, x_j))$ for $i, j = 1, \ldots, n$ as defined in (4.5) with the $r$ suppressed to save notation.

By imposing a second stage prior as in Section 4.3 (see also [1]) on $\psi_0$, a hierarchical Bayes estimator of $\psi$ can be written as

$$\psi_{\mathrm{hb}} = \Xi \Upsilon' H \mathbb{E}[(vI_n + D)^{-1}] H' y, \tag{5.1}$$

where $H$ and $D$ are such that $\Upsilon \Xi \Upsilon' = HDH'$, $H$ is an orthogonal matrix and $D$ is a diagonal matrix.

Hence, using Corollary 2 of [1], a hierarchical Bayes estimator of $f$ is

$$f_{\mathrm{hb}}(x) = (\phi(x)', q(x)')\psi_{\mathrm{hb}}, \tag{5.2}$$

where $\phi(x)$ and $q(x)$ are defined in Section 3.1, $x \in \mathbb{M}$.

5.2. *Hierarchical Bayes estimator as a shrinkage estimator.* Let $\psi_{\text{ls}}$ be the least squares estimator, that is, the solution of the normal equation $\Upsilon'\Upsilon\psi = \Upsilon'y$. We have the following, which shows that the hierarchical Bayes estimator (5.1) is a shrinkage estimator of the least squares solution.

LEMMA 5.1.
$$\psi_{\text{hb}} = (I - \mathbb{E}[v(v\Xi^{-1} + \Upsilon'\Upsilon)^{-1}])\psi_{\text{ls}}.$$

Substituting Lemma 5.1 into (5.2), we have
$$f_{\text{hb}}(x) = (\phi(x)', q(x)')(I - \mathbb{E}[v(v\Xi^{-1} + \Upsilon'\Upsilon)^{-1}])\psi_{\text{ls}},$$

for $x \in \mathbb{M}$.

5.3. *Splines as limits of hierarchical Bayes estimators.* In this section let us compare the hierarchical Bayes estimator (5.2) with the spline estimator of Theorem 3.1. Let us begin by writing
$$\psi_{\text{hb}}(v) = (d'_{\text{hb}}, c'_{\text{hb}}),$$

where $d_{\text{hb}} = \Gamma^{1/2}\Phi'_\Gamma(Q^\Omega_{n,v} + \Phi_\Gamma\Phi'_\Gamma)^{-1}y$, $c_{\text{hb}} = (Q^\Omega_{n,v} + \Phi_\Gamma\Phi'_\Gamma)^{-1}y$, $\Phi_\Gamma = \Phi\Gamma^{1/2}$ and $Q^\Omega_{n,v} = vI + Q^\Omega_n$. Thus, we can rewrite (5.2) as
$$f_{\text{hb}}(x) = \phi(x)'d_{\text{hb}} + q(x)'c_{\text{hb}}$$

for $x \in \mathbb{M}$.

The comparison with the spline estimator of Theorem 3.1 comes from setting $v = n\xi$ and diffusing the parameter $\Gamma$. We use the notation $\|\cdot\|_{\text{op}}$ to denote the usual operator norm. We have the following result.

THEOREM 5.2. *Suppose $x_1, \ldots, x_n \in \mathbb{M}$ locally satisfies the Cox assumptions. If $v = n\xi$, $\xi \asymp n^{-2s/(2s+\dim\mathbb{M})}$ and $\|\Gamma\|_{\text{op}}^{-1} \to 0$, then*
$$E\|f_{\text{hb}} - f\|^2 \ll n^{-2s/(2s+\dim\mathbb{M})}$$

*as $n \to \infty$, where $f \in H_s(\mathbb{M}, M)$, $M > 0$ and $s > (\dim\mathbb{M})(5\dim\mathbb{M} - 2)/4$.*

**6. Application to long period cometary orbits.** Let us illustrate the procedure in the case of the 2-sphere, $S^2$. The data considered in this application consist of directed unit normals of the 658 single-apparition long period cometary orbits found in the catalogue of [33]. The object of interest is the distribution of the directed normals on $S^2$.

This is a well-known directional data set and has been previously analyzed in various ways by Jupp and Spurr [23], Watson [47], Fisher, Lewis and Embleton [15], Wiegert and Tremaine [49] and Mardia and Jupp [32].

Recently, a thorough data analysis is performed in [22], where the reader can obtain more background about this data set. The main conclusions reached in [22] are the following: the data is $SO(2)$-invariant around the North Pole, $(0, 0, 1)'$, or, longitudinally invariant; and, that the data is observed with considerable selection bias whose direct impact is the rejection of spherical uniformity. Indeed, astronomers believe that the intrinsic distribution of the unit directed normals of the long period cometary orbit is spherically uniform and until the Jupp, Kim, Koo and Wiegert [22] paper, it was never really understood why standard directional statistical tests rejected spherical uniformity. It was found that the data has considerable selection bias and when that selection bias is accounted for, Jupp, Kim, Koo and Wiegert [22] show that one can no longer statistically reject null spherical uniformity. In addition to the above findings, with the techniques developed in this paper, we are now able to estimate a longitudinally invariant adaptive estimator of the density of the unit normal vectors to the cometary orbits, that is, the probability density of the observed data with selection bias.

We will do so by using histosplines on $S^2$, since these allow one to calculate densities as a regression problem; see [13]. This is done in the following way. Let

$$\mathcal{S}_{j_1 j_2} = \left[\frac{\pi(j_1 - 1)}{m}, \frac{\pi j_1}{m}\right) \times \left[\frac{2\pi(j_2 - 1)}{m}, \frac{2\pi j_2}{m}\right)$$

for $j_1, j_2 = 1, \ldots, m$. This partitions $S^2$ (with the exception of the South Pole) using (2.3). Thus, we seek the solution to $u \in H_s(S^2)$ that minimizes

$$(6.1) \qquad \sum_{j_1=1}^{m} \sum_{j_2=1}^{m} (y_{j_1 j_2} - L_{j_1,j_2} u)^2 + \xi \int_{S^2} |\Delta^{s/2} u(\omega)|^2 \, d\omega,$$

where $y_{j_1 j_2}$ is the relative frequency of the data in $\mathcal{S}_{j_1 j_2}$ and

$$L_{j_1,j_2} u = \frac{1}{\text{vol}\, S_{j_1 j_2}} \int_{\mathcal{S}_{j_1 j_2}} u(\omega) \, d\omega,$$

for $j_1, j_2 = 1, \ldots, m$. Other than some very minor modifications, the theory would go through exactly as presented in the paper. This is line with what may be called the general spline smoothing problem; see [46], page 10.

The solution to the above general minimization problem (6.1) is

$$u_\xi(x) = \phi(x)' \breve{d} + \breve{q}(x)' \breve{c},$$

where

$$\breve{\Phi} = (L_{j_1,j_2} Y_{kq}),$$

$$\phi(x) = (Y_{kq}(x)),$$

$$\breve{Q}^\iota_{m^2\xi} = [L_{i_1 i_2} L_{j_1 j_2} Q^\iota] + m^2 \xi I_{m^2},$$

$$\breve{d} = (\breve{\Phi}'[\breve{Q}^\iota_{m^2,\xi}]^{-1}\breve{\Phi})^{-1}\breve{\Phi}'[\breve{Q}^\iota_{m^2,\xi}]^{-1}y,$$

$$\breve{c} = [\breve{Q}^\iota_{m^2,\xi}]^{-1}(I_{m^2} - \breve{\Phi}(\breve{\Phi}'[\breve{Q}^\iota_{m^2,\xi}]^{-1}\breve{\Phi})^{-1}\breve{\Phi}'[\breve{Q}^\iota_{m^2,\xi}]^{-1})y,$$

$$Q^\iota(\omega_1,\omega_2) = \sum_{k=K+1}^{\infty} \iota_k \sum_{q=-k}^{k} Y_{kq}(\omega_1)Y_{kq}(\omega_2),$$

$i_1, i_2, j_1, j_2 = 1, \ldots, m$, $|q| \leq k$, $k = 0, 1, \ldots, K$, $x, \omega_1, \omega_2 \in S^2$. Furthermore, $\iota_k = [(k + 1/2)(k + 1)(k + 2)(k + 3)]^{1/2}$ replaces $\lambda_k = k(k + 1)$, which are asymptotically equivalent as $k \to \infty$. This allows us to compute $Q^\iota(\omega_1, \omega_2)$ in closed form, see (4.11) and (4.12), and is the usual way one approaches splines on $S^2$ (cf. [45] and [46], Chapter 2 and [31]).

The above represents the spline solution without any adjustment for $SO(2)$-invariance. If we want to invoke $SO(2)$-invariance together with the hierarchical Bayesian structure, we would need to use the previous $SO(2)$-invariance formulation of Examples 4.2 and 4.4. Thus, all of our parameters are established and, therefore, we can employ the $SO(2)$-invariant adaptive hierarchical Bayes estimator (4.15).

6.1. *Numerical results.* We compute the Bayes estimator, (4.15), for the comet data for several values of $p \in (0, 1)$, $K = 1, 2, \ldots, 10$ and $m = 25$. The Bayes factor (4.18) and Schwarz's criterion (4.19) are given in Table 1 along with the "best" value of $p$ that maximizes the Bayes factor for each value of $K$. We notice that the Bayes factor $\mathfrak{B}_K$ is highest at $K = 6$, while Schwarz's criterion $\mathfrak{S}_K$ is highest at $K = 4$. The two models are very similar and below we compute the Bayes estimator with $K = 6$ and $p = 0.995$. The update value, (4.9), is $p^* = 1 - 1 \times 10^{-10}$, consequently it can be assumed that $p^* = 1$, that is, a posteriori the Bayes model puts all its weight on the $SO(2)$-invariance model and therefore has adapted to the invariance. An explanation of $SO(2)$-invariance for the comet data is given in [22].

In Figure 1 we provide perspective and contour plots of the hierarchical Bayes estimator (with $p^* = 1$ and $K = 6$). The top two panels are the perspective plots, while the bottom two panels are the contour plots. In both sets of plots the domain is taken to be the equal area projection of $S^2$ viewed from the "North Pole" (left panels) and the "South Pole" (right panels). In particular, for each point $(\varphi, \vartheta) \in S^2$, where $\varphi \in [0, 2\pi)$ and $\vartheta \in [0, \pi)$, Lambert's equal area projection is defined by

$$(x,y) = 2\sin\frac{\vartheta}{2}(\cos\varphi, \sin\varphi).$$

By varying $\vartheta \in [0, \pi)$, the 2-sphere is mapped onto a disk of radius 2 with the origin being the "North Pole." The northern hemisphere is thus mapped onto the disk of radius $\sqrt{2}$, while the southern hemisphere is mapped onto an annulus of inner radius $\sqrt{2}$ and outer radius 2. In Figure 1 the left panels take the domain where $\vartheta \in [0, \pi/2]$, so that we are viewing the northern hemisphere of the 2-sphere. We can also reverse the figures so that the South Pole is at the origin with the roles of the southern and northern hemispheres reversed by using

$$(x,y) = 2\sin\frac{(\pi-\vartheta)}{2}(\cos\varphi, \sin\varphi).$$

Indeed, the right panels of Figure 1 are such that one is viewing the "South Pole" at the origin along with just the southern hemisphere. We note that both of these projections preserve area in the sense that $dx\,dy = \sin\vartheta\,d\varphi\,d\vartheta$. In [22] a thorough discussion is given with regard to the equal area projection for the comet data.

As a comparison, in Figure 2 we also plot in a similar way the spline solution where the value of $\xi$ has been chosen by cross-validation; see Remark 3.2.

From Figure 1 one can see that the hierarchical Bayes estimator (4.15) adapts to the $SO(2)$-invariance. In fact, since for $K=6$, $p^* \approx 1$, we have that $\widetilde{\gamma} \approx \widetilde{\gamma}^0$ given by (4.13). Furthermore, the highest concentration of the data is at the North and South Poles which is indicative of the cometary orbits having their orbital planes near the ecliptic plane. The $SO(2)$-invariance is represented by the circular contours. This then represents the estimated probability density of the unit normal vectors of the cometary orbits in the presence of selection bias as explained in [22]. We note that the spline method (see Figure 2) without any adjustment is picking up the peaks at the North and South Poles; however, $SO(2)$-invariance is not particularly distinguishable since the contours do not appear very circular.

Table 1
*Choice of K*

| $K$ | $p$   | $\mathfrak{B}_K$ | $\mathfrak{S}_K$ |
|-----|-------|---------|----------|
| 1   | 0.8   | 0.0000  | 170.8492 |
| 2   | 0.9   | 0.0000  | 232.3569 |
| 3   | 0.8   | 0.0000  | 239.0247 |
| 4   | 0.995 | 49.4024 | 324.4495 |
| 5   | 0.995 | 6.0496  | 278.4142 |
| 6   | 0.995 | 90.0171 | 234.6898 |
| 7   | 0.995 | 7.3810  | 156.4761 |
| 8   | 0.995 | 3.3201  | 143.0733 |
| 9   | 0.995 | 3.0042  | 75.3813  |

**7. Proofs.** In this section we provide the proofs to all of the results in this paper.

7.1. *Proof of Theorem* 3.1. The proof to Theorem 3.1 essentially follows from exhibiting linear independence. For this we need the following result.

LEMMA 7.1. *On a compact connected Riemannian manifold* $\mathbb{M}$ *let*
$$Q(x_i, x) = \phi(x_i)'\phi(x) + Q^\lambda(x_i, x),$$

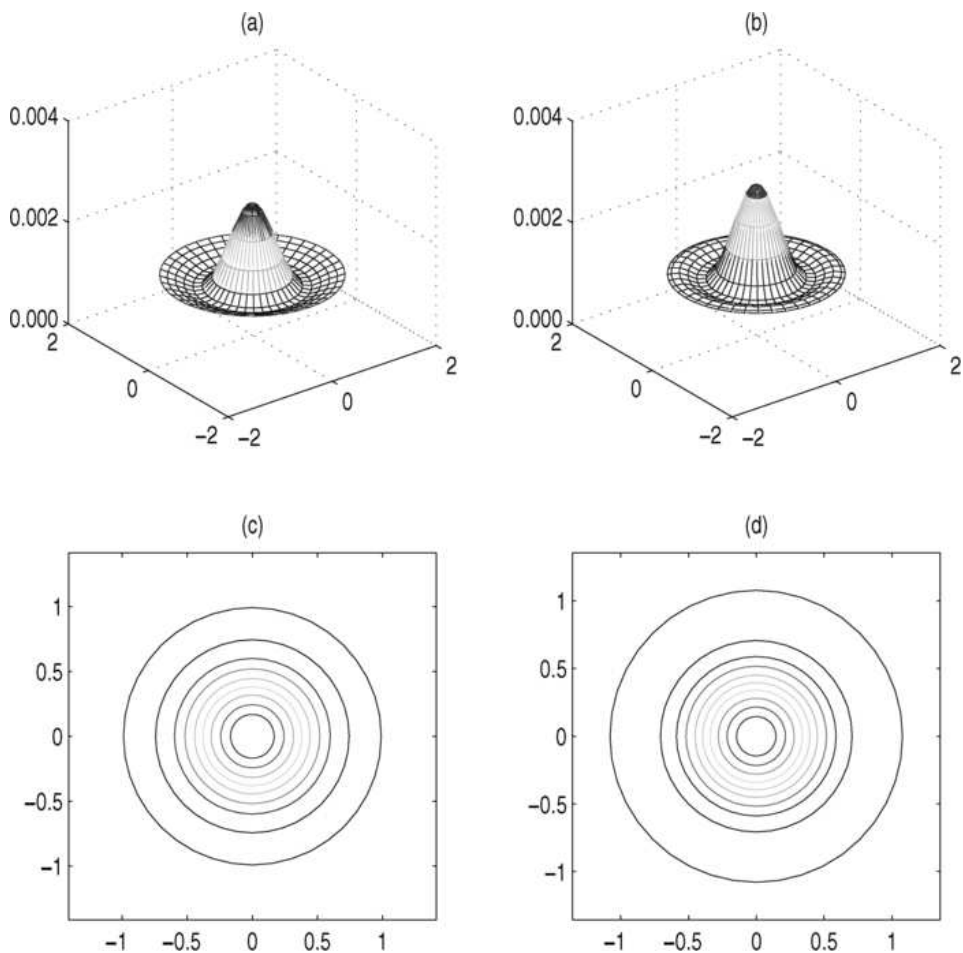

FIG. 1. *Perspective and contour plots of the hierarchical Bayes estimator. The left panels, (a), (c), have the "North Pole" as the origin with just the northern hemisphere as the domain, which is a disk of radius $\sqrt{2}$, while the right panels, (b), (d), have the "South Pole" as the origin with just the southern hemisphere as the domain, which is a disk of radius $\sqrt{2}$.*

where $x_1, \ldots, x_n \in \mathbb{M}$ are distinct and $x \in \mathbb{M}$. Then $\{Q(x_i, x) : i = 1, \ldots, n\}$ is a linearly independent set in $H_s(\mathbb{M})$ for $s > \dim \mathbb{M}$.

PROOF. We first need to regularize the problem as in Lemma 2.3 in [41]. For $p \in \mathbb{M}$, let $(\mathcal{O}_p, \psi_p)$ be a chart; see the Appendix. Now define

$$f_{\epsilon,p}(x) = \begin{cases} \exp\{-\rho(p,x)/[\epsilon - \rho(p,x)]\}, & \text{if } \rho(p,x) \leq \epsilon, \{\rho(p,x) \leq \epsilon\} \subset \mathcal{O}_p, \\ 0, & \text{otherwise,} \end{cases}$$

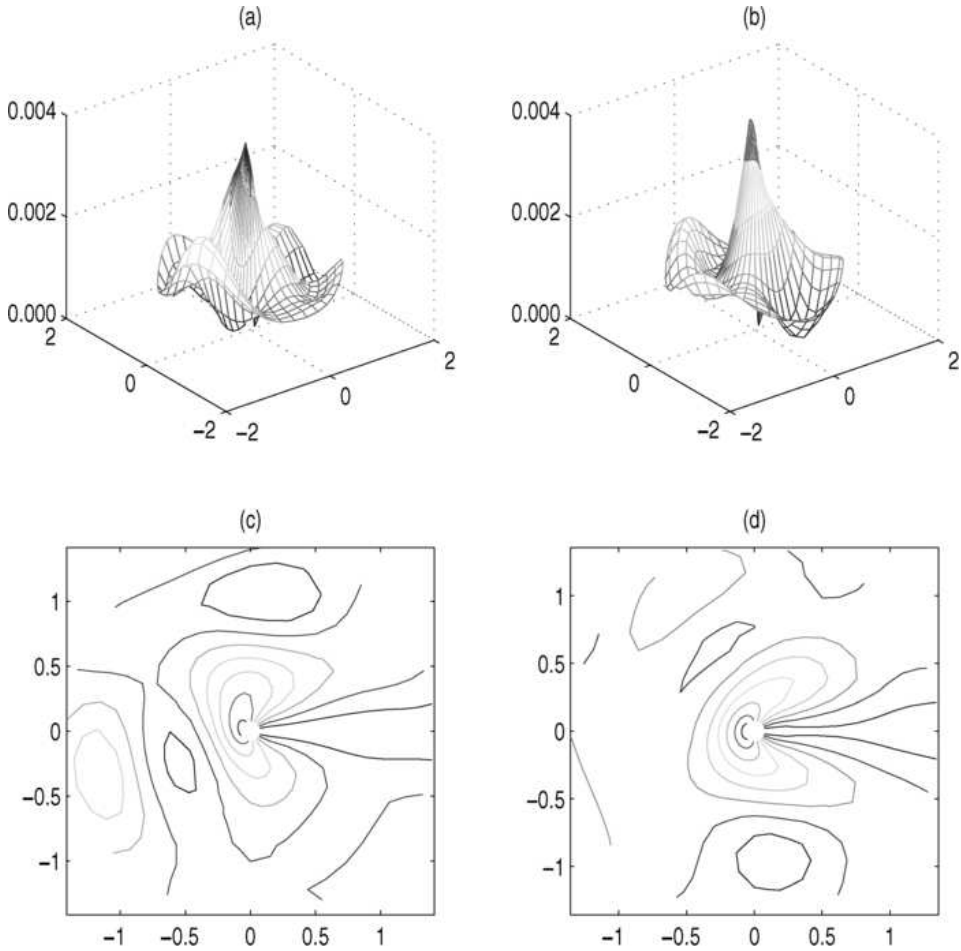

FIG. 2. *Perspective and contour plots of the spline estimator. The left panels,* (a), (c), *have the "North Pole" as the origin with just the northern hemisphere as the domain, which is a disk of radius $\sqrt{2}$, while the right panels,* (b), (d), *have the "South Pole" as the origin with just the southern hemisphere as the domain, which is a disk of radius $\sqrt{2}$.*

where $\rho(\cdot,\cdot)$ is the Riemannian metric; see the Appendix. Notice that we can shrink the compact support of $f_{\epsilon,p}$ around the compact closure of a small open neighborhood around $p \in \mathbb{M}$ just as we would do in the Euclidean case. This will enable us to regularize data.

Consider $x_1, \ldots, x_n$ distinct points in $\mathbb{M}$ and choose $\epsilon$ such that

$$0 < \epsilon < \min_{i_1 \neq i_2} \rho(x_{i_1}, x_{i_2})/2,$$

for $i_1, i_2 = 1, \ldots, n$.

Define

$$u_i(x) = f_{\epsilon, x_i}(x),$$

where $u_i \in C^\infty(\mathbb{M})$ and $u_i \in \mathcal{L}^2(\mathbb{M})$ for $i = 1, \ldots, n$ and $x \in M$. We note that $u_i(x_j) = \delta_{ij}$, where $\delta_{ij}$ denotes the Kronecker delta for $i, j = 1, \ldots, n$.

Suppose that

$$(7.1) \qquad \sum_{i=1}^n \alpha_i Q(x_i, \cdot) = 0.$$

For $g, h \in H_s(\mathbb{M})$ let

$$\langle g, h \rangle = \sum_{k=0}^K \langle \hat{g}_k, \hat{h}_k \rangle_k + \sum_{k=K+1}^\infty \lambda_k^{-s} \langle \hat{g}_k, \hat{h}_k \rangle_k$$

for $s > \dim(\mathbb{M})/2$ and $K > 0$. We note that

$$(7.2) \qquad \langle Q(x_j, \cdot), u_i(\cdot) \rangle = u_i(x_j)$$

for $i, j = 1, \ldots, n$. By applying (7.2) to (7.1), we get

$$0 = \sum_{j=1}^n \alpha_j \langle Q(x_j, \cdot), u_i(\cdot) \rangle = \sum_{j=1}^n \alpha_j u_i(x_j) = \alpha_i$$

for $i = 1, \ldots, n$. Thus, the lemma follows. $\square$

PROOF OF THEOREM 3.1. We note that the $n \times n$ matrix $Q_n^\lambda = (Q^\lambda \times (x_{i_1}, x_{i_2}))$ is positive definite and invertible. Thus, the $n \times n$ matrix $Q_{n,\xi}^\lambda$ is invertible for all $\xi \geq 0$. Now $\Phi$ has rank $\kappa \leq n$ and by Lemma 7.1, for $\xi > 0$ we can apply Theorem 1.3.1 from [46]. The result is the solution to the smoothing problem. $\square$

7.2. *Proof of Theorem* 3.6. The construction of the proof is to localize the arguments over an atlas, do the calculations locally within each chart by applying Theorem 6.2 of Cox [11], and piece together the final argument by using a partition of unity argument; see the Appendix for the technical terms.

Let $\mathcal{O} = \{(\mathcal{O}_\alpha, \psi_\alpha) : \alpha \in \mathcal{A}\}$ be an atlas and consider $\mathcal{P} = \{\delta_\alpha : \alpha \in \mathcal{A}\}$, a partition of unity subordinate to the open sets of the atlas. For a fixed $\delta_\alpha$, let $n_\alpha$ equal the number of $x_j \in \overline{\operatorname{supp} \delta_\alpha}$, $j = 1, \ldots, n$, where the overline in this case means closure. Thus, let $z_{j_i}^\alpha = \psi_\alpha(x_{j_i})$ for $i = 1, \ldots, n_\alpha$.

By the local version of Assumption 3, we can, without loss of generality, assume that the design $x_1, \ldots, x_n \in \mathbb{M}$ is a random sample from the uniform distribution on $\mathbb{M}$. This means that $z_{j_i}^\alpha$, $i = 1, \ldots, n_\alpha$, has distribution $|\partial \psi_\alpha^{-1}(z)| \, dz$, where $|\partial \psi_\alpha^{-1}(z)|$ is the Jacobian of the transformation $\psi_\alpha(x) = z$. Since we are assuming that assumptions 2, 3 and 4 of [11] are locally satisfied, hence, are satisfied on the chart in question, it follows from Theorem 6.2 of [11] that

$$\mathbb{E} \int_{\mathcal{O}_\alpha} |\delta_\alpha(f_\xi^n - f)(x)|^2 \, dx$$
$$= \mathbb{E} \int_{\psi_\alpha(\mathcal{O}_\alpha)} |\delta_\alpha(f_\xi^n - f) \circ \psi_\alpha^{-1}(z)|^2 |\partial \psi_\alpha^{-1}(z)| \, dz$$
$$\ll \xi + n_\alpha^{-1} \xi^{-\dim \mathbb{M}/(2s)}.$$

Now the assumption that there exist constants $0 < c \leq |\partial \psi_\alpha^{-1}(z)| \leq C < \infty$, for all $z \in \overline{\psi_\alpha(\operatorname{supp} \delta_\alpha)}$, allows us to assume

$$\sup_t |\Psi_{n_\alpha}(t) - \Psi_\alpha(t)| \asymp n_\alpha^{-1/\dim \mathbb{M}},$$

where the supremum is taken over $\overline{\psi_\alpha(\operatorname{supp} \delta_\alpha)}$, $\Psi_{n_\alpha}$ is the empirical distribution function of $z_{j_i}^\alpha$, $i = 1, \ldots, n_\alpha$, and $\Psi_\alpha(dt) = |\partial \psi_\alpha^{-1}(t)| \, dt$. Thus, as stated in [11], page 810, in order to satisfy assumption 2 we need

$$s > (\dim \mathbb{M})(5 \dim \mathbb{M} - 2)/4$$

in order to be able to choose $\xi \asymp n_\alpha^{-2s/(2s+\dim \mathbb{M})}$ so that, over the chart $(\mathcal{O}_\alpha, \psi_\alpha)$, we get the asymptotic minimax rate of $n_\alpha^{-2s/(2s+\dim \mathbb{M})}$ as $n_\alpha \to \infty$, $\alpha \in \mathcal{A}$.

The final argument is to use the partition of unity argument to piece together the integration over all of $\mathbb{M}$. Indeed,

$$\mathbb{E} \int_\mathbb{M} |f_\xi^n(x) - f(x)|^2 \, dx = \sum_{\alpha \in \mathcal{A}} \mathbb{E} \int_{\mathcal{O}_\alpha} |\delta_\alpha(f_\xi^n - f)(x)|^2 \, dx$$
$$= \sum_{\alpha \in \mathcal{A}} \mathbb{E} \int_{\psi_\alpha(\mathcal{O}_\alpha)} |\delta_\alpha(f_\xi^n - f) \circ \psi_\alpha^{-1}(z)|^2 |\partial \psi_\alpha^{-1}(z)| \, dz$$

$$\ll \sum_{\alpha \in \mathcal{A}} n_\alpha^{-2s/(2s+\dim \mathbb{M})}$$

$$\ll \left(\sum_{\alpha \in \mathcal{A}} n_\alpha\right)^{-2s/(2s+\dim \mathbb{M})}$$

$$\leq n^{-2s/(2s+\dim \mathbb{M})},$$

where we are using the fact that $\mathcal{A}$ can be taken to be finite, because $\mathbb{M}$ is compact and by the fact that $\sum_{\alpha \in \mathcal{A}} n_\alpha \geq n$.

7.3. *Proof of Theorem* 3.7. The proof is argued along the lines used in [46], Chapter 1. From the Gaussian assumption, we note that

$$\mathbb{E}\tilde{f}_\nu(x)y = \nu\tau^2\phi(x) + \tau^2 q(x),$$

where

$$q(x) = (Q^\lambda(x_1, x), \ldots, Q^\lambda(x_n, x))',$$

$x \in \mathbb{M}$. Furthermore,

$$\mathbb{E}yy' = \nu\tau^2 \Phi\Phi' + \tau^2 Q_n^\lambda + \sigma^2 I_n.$$

Setting $\xi = \sigma^2/n\tau^2$, we have

$$\mathbb{E}(\tilde{f}_\nu(x)|y_1,\ldots,y_n) = \phi(x)'\nu\Phi'(\nu\Phi\Phi' + Q_{n,\xi}^\lambda)^{-1}y + q(x)'(\nu\Phi\Phi' + Q_{n,\xi}^\lambda)^{-1}y$$

for $x \in \mathbb{M}$. Note that

$$\nu\Phi'(\nu\Phi\Phi' + Q_{n,\xi}^\lambda)^{-1} \to (\Phi'[Q_{n,\xi}^\lambda]^{-1}\Phi)^{-1}\Phi'[Q_{n,\xi}^\lambda]^{-1}$$

and

$$(\nu\Phi\Phi' + Q_{n,\xi}^\lambda)^{-1} \to [Q_{n,\xi}^\lambda]^{-1}(I_n - \Phi(\Phi'[Q_{n,\xi}^\lambda]^{-1}\Phi)^{-1}\Phi'[Q_{n,\xi}^\lambda]^{-1})$$

as $\nu \to \infty$.

7.4. *Proof of Lemma* 4.3. Now

$$\mathrm{Cov}(\eta_{i_1}^j, \eta_{i_2}^j) = p^2\tau^2 \sum_{k=K+1}^\infty \langle \phi_k(x_{i_1}), \Omega_k^{j0}\phi_k(x_{i_2})\rangle_k^0$$

$$+ (1-p)^2\tau^2 \sum_{k=K+1}^\infty \langle \phi_k(x_{i_1}), \Omega_k^{j1}\phi_k(x_{i_2})\rangle_k^1$$

for $j = 0, 1$. Since $\Omega_k^{jr} = \mathrm{diag}(\nu_{k\ell}^{jr})$ and $\nu_{k\ell}^{jr} \leq \lambda_k^{-s}$, $\ell = 1, \ldots, \dim \mathcal{E}_k^r$, $k = K+1, \ldots$, $j, r = 0, 1$, we have that

$$|\mathrm{Cov}(\eta_{i_1}^j, \eta_{i_2}^j)| \leq \tau^2 \left| \sum_{k=K+1}^{\infty} \lambda_k^{-s} \langle \phi_k(x_{i_1}), \phi_k(x_{i_2}) \rangle_k \right|$$

$$\leq \tau^2 \left( \sum_{k=K+1}^{\infty} \lambda_k^{-s} \|\phi_k(x_{i_1})\|_k^2 \right)^{1/2} \left( \sum_{k=K+1}^{\infty} \lambda_k^{-s} \|\phi_k(x_{i_2})\|_k^2 \right)^{1/2}$$

$$\leq \tau^2 \sqrt{Z(x_{i_1}, s) Z(x_{i_2}, s)},$$

where $Z(x, s) = \sum_{k>0} \lambda_k^{-s} \|\phi_k(x)\|_k^2$ is the zeta function of $\Delta$. It is known that $Z(x, s)$ is a continuous function of $x$ for fixed $s > \dim(\mathbb{M})/2$ [34]. Since $\mathbb{M}$ is compact, there exists a constant $C(\mathbb{M}, s) < \infty$ depending only on $\mathbb{M}$ and $s$ such that $\sup_{x \in \mathbb{M}} Z(x, s) \leq C(\mathbb{M}, s)$. Hence, $(Q_n^\alpha)_{i,j} \leq C(\mathbb{M}, s)$ for all $i, j = 1, 2, \ldots, n$.

7.5. *Proof of Lemma* 5.1. Using a standard matrix identity (cf. [39], page 151) and omitting the expectation to ease notation, we can write (5.1) as

$$\psi_{\mathrm{hb}} = \Xi \Upsilon' H [(vI_n + D)^{-1}] H' y$$

$$= \Xi \Upsilon' (v H I_n H' + H D H')^{-1} y$$

$$= \Xi \Upsilon' (v I_n + \Upsilon \Xi \Upsilon')^{-1} y$$

$$= [v^{-1} \Xi] \Upsilon (I + \Upsilon [v^{-1} \Xi] \Upsilon')^{-1} y$$

$$= (v \Xi^{-1} + \Upsilon' \Upsilon)^{-1} \Upsilon' y$$

$$= (v \Xi^{-1} + \Upsilon' \Upsilon)^{-1} (\Upsilon' \Upsilon) \psi_{\mathrm{ls}}$$

$$= (v \Xi^{-1} + \Upsilon' \Upsilon)^{-1} ([v \Xi^{-1} + \Upsilon' \Upsilon] - v \Xi^{-1}) \psi_{\mathrm{ls}}$$

$$= (I - v(v \Xi^{-1} + \Upsilon' \Upsilon)^{-1}) \psi_{\mathrm{ls}}.$$

7.6. *Proof of Theorem* 5.2. Again, we will omit the expectation to ease notation. Now

$$\psi_{\mathrm{hb}}(v) = \Xi \Upsilon' H (vI_n + D)^{-1} H' y$$

$$= \Xi \Upsilon' (v H I_n H' + H D H')^{-1} y$$

$$= \Xi \Upsilon' (v I_n + \Upsilon \Xi \Upsilon')^{-1} y$$

$$= \begin{pmatrix} \Gamma & 0 \\ 0 & Q_n^\Omega \end{pmatrix} \begin{pmatrix} \Phi' \\ I \end{pmatrix} (vI_n + Q_n^\Omega + \Phi\Gamma\Phi')^{-1} y$$

$$= \begin{pmatrix} \Gamma\Phi' \\ Q_n^\Omega \end{pmatrix} (vI + Q_n^\Omega + \Phi\Gamma^{1/2}\Gamma^{1/2}\Phi')^{-1} y$$

$$= \begin{pmatrix} \Gamma\Phi' \\ Q_n^\Omega \end{pmatrix} (Q_{n,v}^\Omega + \Phi_\Gamma\Phi_\Gamma')^{-1} y,$$

where $Q_v^\Omega = vI + Q_n^\Omega$ and $\Phi_\Gamma = \Phi\Gamma^{1/2}$. Hence, the first $\kappa$ components of $\psi_{\text{hb}}(v)$ are given by

$$\Gamma\Phi'(Q_{n,v}^\Omega + \Phi_\Gamma\Phi_\Gamma')^{-1} y = \Gamma^{1/2}[\Phi\Gamma^{1/2}]'(Q_{n,v}^\Omega + \Phi_\Gamma\Phi_\Gamma')^{-1} y$$
$$= \Gamma^{1/2}\Phi_\Gamma'(Q_{n,v}^\Omega + \Phi_\Gamma\Phi_\Gamma')^{-1} y,$$

and the last $n$ ones are

$$Q_{n,v}^\Omega(Q_{n,v}^\Omega + \Phi_\Gamma\Phi_\Gamma')^{-1} y.$$

Using Angers and Delampady [1], the hierarchical Bayesian equivalent of $d$ of Theorem 3.1 is then

$$\begin{aligned}
d_{\text{hb}} &= \Gamma^{1/2}\Phi_\Gamma'(Q_{n,v}^\Omega + \Phi_\Gamma\Phi_\Gamma')^{-1} y \\
&= \Gamma^{1/2}[(Q_{n,v}^\Omega + \Phi_\Gamma\Phi_\Gamma')^{-1}\Phi_\Gamma]' y \\
&= \Gamma^{1/2}[(I + \Phi_\Gamma'(Q_{n,v}^\Omega)^{-1}\Phi_\Gamma)\Phi_\Gamma'(Q_{n,v}^\Omega)^{-1}] y \\
&= \Gamma^{1/2}(I + \Gamma^{1/2}\Phi'\Phi\Gamma^{1/2})^{-1}\Gamma^{1/2}\Phi'(Q_{n,v}^\Omega)^{-1} y \\
&= (\Gamma^{-1} + \Phi'(Q_{n,v}^\Omega)^{-1}\Phi)^{-1}\Phi'(Q_{n,v}^\Omega)^{-1} y \\
&\to d
\end{aligned}$$

if $\|\Gamma\|_{\text{op}}^{-1} \to 0$ and $v = n\xi$. Similarly, the hierarchical Bayesian equivalent of $c$ of Theorem 3.1 is given by

$$\begin{aligned}
c_{\text{hb}} &= [Q_n^\Omega]^{-1} Q_n^\Omega (Q_{n,v}^\Omega + \Phi_\Gamma\Phi_\Gamma')^{-1} y \\
&= (Q_{n,v}^\Omega + \Phi_\Gamma\Phi_\Gamma')^{-1} y \\
&= (Q_{n,v}^\Omega)^{-1}(I + \Phi_\Gamma\Phi_\Gamma'(Q_{n,v}^\Omega)^{-1})^{-1} y \\
&= (Q_{n,v}^\Omega)^{-1}[I - \Phi_\Gamma\Phi_\Gamma'(I + (Q_{n,v}^\Omega)^{-1}\Phi_\Gamma\Phi_\Gamma')^{-1}(Q_{n,v}^\Omega)^{-1}] y \\
&= (Q_{n,v}^\Omega)^{-1}[I - \Phi_\Gamma\Phi_\Gamma'(Q_{n,v}^\Omega + \Phi_\Gamma\Phi_\Gamma')^{-1} Q_{n,v}^\Omega (Q_{n,v}^\Omega)^{-1}] y \\
&= (Q_{n,v}^\Omega)^{-1}[I - \Phi\Gamma^{1/2}\Phi_\Gamma'(Q_{n,v}^\Omega + \Phi_\Gamma\Phi_\Gamma')^{-1}] y \\
&= (Q_{n,v}^\Omega)^{-1}[I - \Phi(\Gamma^{-1} + \Phi'(Q_{n,v}^\Omega)^{-1}\Phi)^{-1}\Phi'(Q_{n,v}^\Omega)^{-1}] y \\
&\to c
\end{aligned}$$

if $\|\Gamma\|_{\text{op}}^{-1} \to 0$ and $v = n\xi$.

Consequently, the spline estimator of Theorem 3.1 is a limiting case of the hierarchical Bayes estimator (5.2)

$$f_{\text{hb}}(x) \to f_\xi^n(x), \tag{7.3}$$

for $x \in \mathbb{M}$, when $v = n\xi$ as $\|\Gamma\|_{\text{op}}^{-1} \to 0$. Now we know that

$$E\|f_{\text{hb}} - f\| \leq E\|f_{\text{hb}} - f_\xi^n\| + E\|f_\xi^n - f\|.$$

By (7.3), we know that $E\|f_{\text{hb}} - f_\xi^n\| \to 0$ as $\|\Gamma\|_{\text{op}}^{-1} \to 0$ for $v = n\xi$. Thus, the result follows.

## APPENDIX

In this appendix we gather together some further technical discussions that concern $\mathbb{M}$ and that are directly or indirectly used in the paper. Consider a smooth curve $\gamma:[a,b] \to \mathbb{M}$, with $a < b$, and let $\gamma'(t)$ denote its first derivative for some $t \in (a,b)$. Then the length of $\gamma$ is defined through the Riemannian structure as

$$l(\gamma) = \int_a^b g_{\gamma(t)}(\gamma'(t), \gamma'(t))^{1/2} \, dt.$$

Since we are assuming that $\mathbb{M}$ is connected, hence, for any two points $p, q \in \mathbb{M}$, we can find a curve in $\mathbb{M}$ that joins them in $\mathbb{M}$, we can define a metric on $\mathbb{M}$ by

$$\rho(p, q) = \inf\{l(\gamma) : \gamma \text{ joining } p \text{ and } q\},$$

$p, q \in \mathbb{M}$. This metric is called the Riemannian distance or metric which makes $(\mathbb{M}, \rho)$ a metric space ([7], page 39).

For $p \in \mathbb{M}$, let $(\mathcal{O}_\alpha, \psi_\alpha)$ be a chart, that is, $\mathcal{O}_\alpha \subset \mathbb{M}$ is an open set with $p \in \mathcal{O}_\alpha$ and $\psi_\alpha : \mathcal{O}_\alpha \to \psi_\alpha(\mathcal{O}_\alpha) \subset \mathbb{R}^{\dim \mathbb{M}}$ is a diffeomorphism. A collection of charts $\{(\mathcal{O}_\alpha, \psi_\alpha) : \alpha \in \mathcal{A}\}$ is called an atlas if $\bigcup_\alpha \mathcal{O}_\alpha = \mathbb{M}$. Since $\mathbb{M}$ is assumed to be compact, we can take $\mathcal{A}$ to be a finite set. Thus, the open sets of an atlas form a finite open cover of $\mathbb{M}$. Consider a collection $\mathcal{P}$ of nonnegative functions $\delta_\alpha$ whose support is contained in $\mathcal{O}_\alpha$ and that has the property that $1 = \sum_\alpha \delta_\alpha$. We call such a collection a partition of unity subordinate to the open cover, and for some integrable functions $f, g : \mathbb{M} \to \mathbb{R}$, integration is defined by

$$\int_{\mathbb{M}} f(x) g(x) \, dx = \sum_{\alpha \in \mathcal{A}} \int_{\psi_\alpha(\mathcal{O}_\alpha)} f \circ \psi_\alpha^{-1}(z) (\delta_\alpha g) \circ \psi_\alpha^{-1}(z) \, dz. \tag{A.1}$$

We note that this definition is well defined, hence, is independent of the choice of open cover ([7], page 6).

A famous formula due to Hermann Weyl states

$$\lim_{k\to\infty} \frac{\lambda_k}{(\sum_{j=0}^{k} \dim \mathcal{E}_j)^{2/\dim \mathbb{M}}} = \mathcal{W}^{-2/\dim \mathbb{M}}$$

where

$$\text{(A.2)} \qquad \mathcal{W} = \frac{\text{vol}\,\mathbb{M}}{(2\sqrt{\pi})^{\dim \mathbb{M}} \Gamma(1 + \dim \mathbb{M}/2)}$$

([7], page 9). We note that the $\mathcal{W}$ appearing above is the same quantity which appears in the minimax constant of Theorem 3.4.

**Acknowledgments.** We would like to thank the Editor, the Associate Editor and two anonymous referees for their careful reviews of this manuscript. We are especially grateful to one of the referees who provided very insightful comments with respect to the geometry and statistics of this paper.

## REFERENCES

[1] ANGERS, J.-F. and DELAMPADY, M. (1992). Hierarchical Bayesian curve fitting and smoothing. *Canad. J. Statist.* **20** 35–49. MR1173059

[2] ANGERS, J.-F. and DELAMPADY, M. (1997). Hierarchical Bayesian curve fitting and model choice for spatial data. *Sankhyā Ser. B* **59** 28–43. MR1733376

[3] ANGERS, J.-F. and DELAMPADY, M. (2001). Bayesian nonparametric regression using wavelets. *Sankhyā Ser. A* **63** 287–308. MR1897044

[4] BERAN, R. J. (1968). Testing for uniformity on a compact homogeneous space. *J. Appl. Probability* **5** 177–195. MR0228098

[5] BERAN, R. (1979). Exponential models for directional data. *Ann. Statist.* **7** 1162–1178. MR0550142

[6] BERGER, J. O. (1985). *Statistical Decision Theory and Bayesian Analysis*, 2nd ed. Springer, New York. MR0804611

[7] CHAVEL, I. (1984). *Eigenvalues in Riemannian Geometry*. Academic Press, Orlando, FL. MR0768584

[8] CHIKUSE, Y. (2003). *Statistics on Special Manifolds. Lecture Notes in Statist.* **174**. Springer, New York. MR1960435

[9] CHIKUSE, Y. and JUPP, P. E. (2004). A test of uniformity on shape spaces. *J. Multivariate Anal.* **88** 163–176. MR2021868

[10] CHIRIKJIAN, G. S. and KYATKIN, A. B. (2001). *Engineering Applications of Noncommutative Harmonic Analysis*. CRC Press, Boca Raton, FL. MR1885369

[11] COX, D. D. (1984). Multivariate smoothing spline functions. *SIAM J. Numer. Anal.* **21** 789–813. MR0749371

[12] COX, D. D. (1988). Approximation of method of regularization estimators. *Ann. Statist.* **16** 694–712. MR0947571

[13] DYN, N. and WAHBA, G. (1982). On the estimation of functions of several variables from aggregated data. *SIAM J. Math. Anal.* **13** 134–152. MR0641546

[14] EFROMOVICH, S. (2000). On sharp adaptive estimation of multivariate curves. *Math. Methods Statist.* **9** 117–139. MR1780750

[15] FISHER, N. I., LEWIS, T. and EMBLETON, B. J. J. (1993). *Statistical Analysis of Spherical Data*. Cambridge Univ. Press. MR1247695

[16] GINÉ, E. M. (1975). Invariant tests for uniformity on compact Riemannian manifolds based on Sobolev norms. *Ann. Statist.* **3** 1243–1266. MR0388663

[17] HANNA, M. S. and CHANG, T. (2000). Fitting smooth histories to rotation data. *J. Multivariate Anal.* **75** 47–61. MR1787401

[18] HEALY, D. M., HENDRIKS, H. and KIM, P. T. (1998). Spherical deconvolution. *J. Multivariate Anal.* **67** 1–22. MR1659108

[19] HEALY, D. M. and KIM, P. T. (1996). An empirical Bayes approach to directional data and efficient computation on the sphere. *Ann. Statist.* **24** 232–254. MR1389889

[20] HENDRIKS, H. (1990). Nonparametric estimation of a probability density on a Riemannian manifold using Fourier expansions. *Ann. Statist.* **18** 832–849. MR1056339

[21] JUPP, P. E. and KENT, J. T. (1987). Fitting smooth paths to spherical data. *Appl. Statist.* **36** 34–46. MR0887825

[22] JUPP, P. E., KIM, P. T., KOO, J.-Y. and WIEGERT, P. (2003). The intrinsic distribution and selection bias of long-period cometary orbits. *J. Amer. Statist. Assoc.* **98** 515–521. MR2011668

[23] JUPP, P. E. and SPURR, B. D. (1983). Sobolev tests for symmetry of directional data. *Ann. Statist.* **11** 1225–1231. MR0720267

[24] KIM, P. T. (1998). Deconvolution density estimation on $SO(N)$. *Ann. Statist.* **26** 1083–1102. MR1635446

[25] KIM, P. T. and KOO, J.-Y. (2000). Directional mixture models and optimal estimation of the mixing density. *Canad. J. Statist.* **28** 383–398. MR1792056

[26] KIM, P. T. and KOO, J.-Y. (2002). Optimal spherical deconvolution. *J. Multivariate Anal.* **80** 21–42. MR1889831

[27] KIM, P. T. and RICHARDS, D. ST. P. (2001). Deconvolution density estimation on compact Lie groups. In *Algebraic Methods in Statistics and Probability* (M. A. G. Viana and D. St. P. Richards, eds.) 155–171. Amer. Math. Soc., Providence, RI. MR1873674

[28] KLEMELÄ, J. (1999). Asymptotic minimax risk for the white noise model on the sphere. *Scand. J. Statist.* **26** 465–473. MR1712031

[29] LEE, J. and RUYMGAART, F. H. (1996). Nonparametric curve estimation on Stiefel manifolds. *J. Nonparametr. Statist.* **6** 57–68. MR1382072

[30] LINDLEY, D. V. and SMITH, A. F. M. (1972). Bayes estimates for the linear model (with discussion). *J. Roy. Statist. Soc. Ser. B* **34** 1–41. MR0415861

[31] LUO, Z. (1998). Backfitting in smoothing spline ANOVA. *Ann. Statist.* **26** 1733–1759. MR1673276

[32] MARDIA, K. V. and JUPP, P. E. (2000). *Directional Statistics*. Wiley, Chichester. MR1828667

[33] MARSDEN, B. G. and WILLIAMS, G. V. (1993). *Catalogue of Cometary Orbits*, 8th ed. Minor Planet Center, Smithsonian Astrophysical Observatory, Cambridge, MA.

[34] MINAKSHISUNDARAM, S. and PLEIJEL, Å. (1949). Some properties of the eigenfunctions of the Laplace operator on Riemannian manifolds. *Canadian J. Math.* **1** 242–256. MR0031145

[35] PINSKER, M. S. (1980). Optimal filtering of square integrable signals in Gaussian white noise. *Problems Inform. Transmission* **16** (1) 52–68. MR0624591

[36] PRENTICE, M. J. (1987). Fitting smooth paths to rotation data. *Appl. Statist.* **36** 325–331. MR0918853

[37] ROBERT, C. P. (2001). *The Bayesian Choice*, 2nd ed. Springer, New York. MR1835885

[38] Schwarz, G. (1978). Estimating the dimension of a model. *Ann. Statist.* **6** 461–464. MR0468014
[39] Searle, S. R. (1982). *Matrix Algebra Useful for Statistics*. Wiley, New York. MR0670947
[40] Speckman, P. (1985). Spline smoothing and optimal rates of convergence in nonparametric regression models. *Ann. Statist.* **13** 970–983. MR0803752
[41] Taijeron, H., Gibson, A. and Chandler, C. (1994). Spline interpolation and smoothing on hyperspheres. *SIAM J. Sci. Comput.* **15** 1111–1125. MR1289156
[42] Turski, J. (1998). Harmonic analysis on $SL(2,\mathbb{C})$ and projectively adapted pattern representation. *J. Fourier Anal. Appl.* **4** 67–91. MR1650956
[43] Turski, J. (2000). Projective Fourier analysis for patterns. *Pattern Recognition* **33** 2033–2043.
[44] van Rooij, A. C. M. and Ruymgaart, F. H. (1991). Regularized deconvolution on the circle and the sphere. In *Nonparametric Functional Estimation and Related Topics* (G. Roussas, ed.) 679–690. Kluwer, Dordrecht. MR1154359
[45] Wahba, G. (1981). Spline interpolation and smoothing on the sphere. *SIAM J. Sci. Statist. Comput.* **2** 5–16. [Erratum (1982) **3** 385–386.] MR0618629
[46] Wahba, G. (1990). *Spline Models for Observational Data*. SIAM, Philadelphia. MR1045442
[47] Watson, G. S. (1983). *Statistics on Spheres*. Wiley, New York. MR0709262
[48] Wellner, J. A. (1979). Permutation tests for directional data. *Ann. Statist.* **7** 929–943. MR0536498
[49] Wiegert, P. and Tremaine, S. (1999). The evolution of long-period comets. *Icarus* **137** 84–121.

Département de Mathématiques
et de Statistique
Université de Montréal
Montréal, Québec
Canada H3C 3J7
e-mail: angers@dms.umontreal.ca

Department of Mathematics
and Statistics
University of Guelph
Guelph, Ontario
Canada N1G 2W1
e-mail: pkim@uoguelph.ca